\documentclass [reqno]{amsart}
\usepackage {mathptmx}
\usepackage {amsmath,amssymb}
\usepackage {color,graphicx}

\newcommand{\ignore}[1]{} 

\renewenvironment {enumerate}%
  {\rule{1mm}{0mm}\begin {oldenumerate}%
    \parskip1ex plus0.5ex \itemsep 0mm \parindent 0mm}%
  {\end {oldenumerate}}

\renewenvironment {itemize}%
  {\rule{1mm}{0mm}\begin {olditemize}%
    \parskip1ex plus0.5ex \itemsep 0mm \parindent 0mm}%
  {\end {olditemize}}

\theoremstyle {plain}
\newtheorem {theorem}{Theorem}[section]
\newtheorem {proposition}[theorem]{Proposition}
\newtheorem {lemma}[theorem]{Lemma}
\newtheorem {lemdef}[theorem]{Lemma and Definition}
\newtheorem {corollary}[theorem]{Corollary}
\newtheorem{conj}[theorem]{Question}
\theoremstyle {definition}
\newtheorem {definition}[theorem]{Definition}
\newtheorem {remdef}[theorem]{Definition and Remark}
\newtheorem {remark}[theorem]{Remark}
\newtheorem {construction}[theorem]{Construction}
\newtheorem {convention}[theorem]{Convention}
\newtheorem {example}[theorem]{Example}

\newcommand{\N}{\mathbb{N}}
\newcommand{\Z}{\mathbb{Z}}
\newcommand{\Q}{\mathbb{Q}}
\newcommand{\R}{\mathbb{R}}
\newcommand{\C}{\mathbb{C}}

\newcommand{\PP}{\mathbb{P}}

\newcommand{\M}{\mathcal{M}}

\newcommand{\bet}[1] {|{#1}|}
\newcommand {\mult}{\mathop {\rm mult}\nolimits}

\newcommand{\be}{\begin{equation}}
\newcommand{\ee}{\end{equation}}
\newcommand{\bi}{\begin{itemize}}
\newcommand{\ei}{\end{itemize}}

\newcommand{\lk}{\left\{ }
\newcommand{\rk}{\right\} }

\newcommand{\tn}{\textnormal}

\newcommand{\ev}{\tn{ev}}

\begin {document}

  \title [Irreducible cycles and points in special position]%
  {Irreducible cycles and points in special position in moduli spaces for 
  tropical curves}
\author{Andreas Gathmann and Franziska Schroeter}
\address {Andreas Gathmann, Fachbereich Mathematik, TU Kaiserslautern, Postfach
  3049, 67653 Kaiserslautern, Germany}
\email {andreas@mathematik.uni-kl.de}
\address {Franziska Schroeter, Cluster of Excellence for Computer Science,
  Universit\"at des Saarlandes, Postfach 151150, 66041 Saarbr\"ucken, Germany}
  \email{schroeter@math.uni-sb.de}

\makeatletter
\@namedef{subjclassname@2010}{\textup{2010} Mathematics Subject Classification}
\makeatother


\begin {abstract}
  In the first part of this paper, we discuss the notion of irreducibility of
  cycles in the moduli spaces of $n$-marked rational tropical curves. We prove
  that Psi-classes and vital divisors are irreducible, and that locally
  irreducible divisors are also globally irreducible for $ n \le 6 $. In the
  second part of the paper, we show that the locus of point configurations in $
  (\R^2)^n $ in special position for counting rational plane curves (defined in
  two different ways) can be given the structure a tropical cycle of
  codimension $1$. In addition, we compute explicitly the weights of this
  cycle.
\end {abstract}

\maketitle

  \section {Introduction} \label{l1}

In classical algebraic geometry, a divisor on a nonsingular algebraic variety
$X$ over $\C$ is a (formal) finite $\Z$-linear combination of irreducible
subvarieties of codimension $1$. Interpreting a divisor as an algebraic variety
(with multiplicities attached to its components) it is of course irreducible if
and only if it consists of only one subvariety. In tropical geometry however,
we define a divisor in a tropical fan $X$ to be, roughly speaking, a weighted
subfan of codimension one --- so in the case of non-negative weights it is just
an embedded tropical variety. As there is no unique decomposition of a such a
variety into irreducible components (see Remark \ref{rem-class}), the situation
is much more complicated here. The notion of irreducible tropical fan was
introduced in \cite{GKM} and can easily be generalized to tropical divisors in
$X$ (Definition \ref{def-irred}). A divisor is then irreducible if and only if
its weight function is unique up to a global multiple (Lemma \ref{criterion}).
We explore this notion when the ambient space is $ X=\mathcal{M}_{0,n} $, the
parameter space of tropical rational $n$-marked tropical curves \cite{GKM,M},
or equivalently the tropical Grassmannian $\mathcal{G}(2,n)$ modulo its
lineality space \cite{SS}. We prove for example that tropical Psi-classes and
vital divisors are irreducible. Note that the corresponding statement is also
true in algebraic geometry for the vital divisors \cite{K}, where they form the
main objects in the intersection theory on the moduli spaces
$\overline{M}_{0,n}$ of stable curves. In contrast, Psi-classes in algebraic
geometry cannot be defined as divisors, but only as divisor classes on
$\overline{M}_{0,n}$ \cite{AC}, and hence the notion of irreducibility does not
make sense for them. They are useful for the construction of enumerative
invariants like descendant Gromov-Witten numbers as in \cite{KoM}. Furthermore
we show that for $ n \le 6 $ every divisor in $\mathcal{M}_{0,n}$ that is
locally irreducible is already globally irreducible (Proposition
\ref{prop-m06}).

The second part of this paper, which is almost independent of the first, deals
with the tropical notion of points in special position for enumerative
problems. If one wants to count for instance curves of degree $d$ in $ \PP^2 $
passing through $ n=3d-1 $ given points, these point conditions should be
generic enough to obtain at least a finite result that is independent of the
conditions. We call the set of these generic point configurations in
$(\PP^2)^n$ the set of \textit{points in general position}. There exist several
definitions of this notion --- we can for example just require the number of
curves to be finite, or in addition that the curves are all smooth. A common
feature of these definitions is that they lead to an open dense subset of
$(\PP^2)^n$ of points in general position. The complement of this set in
$(\PP^2)^n$ is called the set of \textit{points in special position}. We
explore this notion in the corresponding tropical context. It turns out that
this set is the support of a tropical fan in $(\R^2)^n$, more precisely of a
push-forward of a divisor in the tropical moduli space corresponding to the
classical moduli space $ \overline{\mathcal{M}}_{0,3d-1} (\PP^2,d)$ of stable
maps (Propositions \ref{v1} and \ref{v2}). Finally, we compute the weights of
the top-dimensional cones in this tropical fan by means of elementary divisors
of certain matrices (Corollary \ref {resum}).

We would like to thank Hannah Markwig for inspiring and helpful discussions.
Part of this work was done at the Mathematical Sciences Research Institute
(MSRI) in Berkeley, CA, USA, during the one-semester program on tropical
geometry in fall 2009. The authors would like to thank the MSRI for
hospitality and support. In particular, Andreas Gathmann was supported by
the Simons Professorship of the MSRI.   

  \section{(Irreducible) tropical subfans and cycles in $\M_{0,n}$} \label{l2}

In this section, we want to verify that the codimension-$k$ skeletons of the
moduli spaces $\M_{0,n}$ of rational $n$-marked abstract tropical curves are
tropical subvarieties, for suitable weight functions. Furthermore, it turns out
that the same holds for codimension-$k$ skeletons of Psi-classes $\psi_i$
living in $\M_{0,n}$. This will be done by checking the balancing condition
explicitly in terms of coordinate vectors of $\M_{0,n}$. Inspired by this
method, we will show in the subsequent part that tropical Psi-classes $\psi_i$
are irreducible. In addition, we will prove that the tropical vital divisors
$D^S$ in $\M_{0,n}$ are also irreducible. As a by-product we will give a new
characterization of irreducible tropical fans and cycles.

Let us start by reviewing the definitions and notions indispensable below and
by fixing the notation. More details about the following notions can be found in
\cite[Definitions 2.4, 2.6. and 2.8]{GKM}.

\begin{definition}[Tropical (sub-)fans, codimension-$k$ skeleton] \label {defo}
   Pick some $r \geq 0$ and a lattice $ \Lambda
  \cong \Z^r $; set $V:=\Lambda \otimes_\Z \R$.
  \begin {enumerate}
  \item Let $X$ be a fan of pure dimension $N$ in $V$. For $ 0 \le k \le N $ we
    denote the collection of its $k$-dimensional cones by $ X^{(k)} $. The
    vector subspace of $V$ spanned by a cone $ \sigma $ of $X$ will be denoted
    $ V_\sigma $. We call $X$ a \textit {tropical fan} if it is equipped with a
    weight function $ w: X^{(N)} \to \N $ such that the balancing condition
      \[ \sum_{\sigma>\tau} w(\sigma) \, v_{\sigma/\tau} = 0 \quad \in V/V_\tau
      \]
    holds for all $ (N-1) $-dimensional cones $ \tau $, where $ v_{\sigma/\tau}
    $ is the primitive normal vector of $ \sigma $ relative to $ \tau $. The
    union of all cones of $X$ will be written as $ |X| \subset V $.
  \item A \textit{(tropical) subfan} $Y$ of a (tropical) fan $X$ in $V$ is
    itself a (tropical) fan $Y$ in $V$ with the property that each cone of $Y$
    is contained in a cone of $X$. Note that the weight function of $Y$ is not
    necessarily inherited by $X$.
  \item The \textit{codimension-$k$ skeleton} $X^{(N-k)}_{\tn{sk}}$ of a fan
    $X$ in $V$ of pure dimension $N$ (for $ 0 \leq k \leq N$) consists of all
    cones of dimension at most $ N-k $ in $X$. It is a pure-dimensional fan of
    dimension $N-k$, however with no canonical weight function associated to
    it.
  \end {enumerate}
\end{definition}

\begin{construction}[Moduli space $\M_{0,n}$ of rational $n$-marked abstract
    tropical curves] \label {moduli}
  A \textit {rational $n$-marked (abstract) tropical curve} is a tuple
  $ (\Gamma;x_1,\ldots,x_n) $, where $ \Gamma $ is a connected metric tree
  graph all of whose vertices have valence at least 3, and where $ x_1,\dots,
  x_n $ is a labeling of its unbounded edges. In the following, we will sometimes
  write marked edge $x_i$ or marked end $x_i$ for the unbounded edge with label $x_i$. An
  \textit{isomorphism} of rational tropical curves $(\Gamma;x_1,\ldots,x_n) $ and
  $(\tilde{\Gamma};\tilde{x}_1,\ldots,\tilde{x}_n)$ is a homeomorphism $\Gamma \to
  \tilde{\Gamma}$ sending $x_i$ to $\tilde{x}_i$ and mapping the edges of
  $\Gamma$ bijectively onto edges of $\tilde{\Gamma}$ by affine maps of slope
  $\pm 1$. As in \cite{GKM}, we will denote by $\M_{0,n}$ the parameter space
  of all such curves modulo isomorphisms. For many purposes it plays the
  tropical role of the classical moduli space of rational $n$-marked stable
  curves --- even if it is not compact.

  The moduli space $ \M_{0,n} $ can be given the structure of a tropical fan
  (see Remark \ref{combtypes} for the fan structure) of dimension $n-3$ in a
  quotient space of $ \R^{\binom{n}{2}}$ (\cite[Theorem 3.7]{GKM});
  in fact, it can be described as the tropical Grassmannian
  $\mathcal{G}(2,n)$ modulo its lineality space (see \cite[Theorem 3.4]{SS} or
\cite[Remark 3.9]{GKM}).
  More precisely, for each $n$-marked tropical curve $ (\Gamma;
  x_1,\dots,x_n) $ and $ 1 \le i < j \le n $ let $ \tn{dist}_\Gamma(x_i,x_j) $
  be the distance between the unbounded edges $x_i$ and $x_j$ in $\Gamma$. We
  thus get a map
  \begin {align*}
    \tilde v: \qquad\qquad \M_{0,n}  &\to     \R^{\binom n2} \\
              (\Gamma;x_1,\dots,x_n) &\mapsto
                (\tn{dist}_\Gamma(x_i,x_j))_{i<j}
  \end {align*}
  where we choose the lexicographic ordering of the pairs $ (i,j) $ for the
  coordinates in $ \R^{\binom n2} $. We will call $ \tilde v(\Gamma;
  x_1,\dots,x_n) $ the \textit {distance vector} of $ (\Gamma;x_1,\dots,x_n) $.

  The following vectors in $ \R^{\binom n2} $ will be of particular importance:
  let $ I \subset \{1,\dots,n\} $ be any subset, and denote by $ \tilde v(I)
  \in \R^{\binom n2} $ the vector whose $ (i,j) $-coordinate is equal to 1
  if $I$ contains exactly one of the numbers $i$ and $j$, and 0 otherwise. Note
  that $ \tilde v(I^c) = \tilde v(I) $, where $ I^c $ denotes the complement of
  $I$ in $ \{1,\dots,n\} $.

  The vectors $ d_i := \tilde v(\{i\}) $ for $ i=1,\dots,n $ form a basis of
  the so-called lineality space mentioned above; i.e.\ by taking the quotient
  by this subspace we obtain a map
  \begin {align*}
    v: \qquad\qquad \M_{0,n}  &\to
         \R^{\binom n2}/\langle d_1,\dots,d_n\rangle \\
       (\Gamma;x_1,\dots,x_n) &\mapsto
         \overline {(\tn{dist}_\Gamma(x_i,x_j))_{i<j}}
  \end {align*}
  that embeds $ \M_{0,n} $ as a tropical fan in $ \R^{\binom n2}/\langle
  d_1,\dots,d_n\rangle $ \cite[Theorem 3.7]{GKM}. For this structure of a tropical
fan the weights of
  all top-dimensional cones are chosen to be $1$, and the underlying lattice
  $\Lambda$ is taken to be the one generated by the classes $ v(I) $ of the
  vectors $ \tilde v(I) $ for all $ I \subset \{1,\dots,n\} $ modulo the
  lineality space. In the remainder, we will always view $ \M_{0,n} $ as a
  tropical fan with this embedding.
\end {construction}

\begin {remark} \label {combtypes}
  The cones of the tropical fan $\M_{0,n}$ correspond to curves of the same
  \textit{combinatorial type}, i.e.\ curves which differ just by the lengths of their
  bounded edges. They are open convex polyhedra of dimension
    $$ n-3- \sum_{\tn{ vertices } V \tn{ in } \Gamma}
       (\tn{val } V-3), $$
  where $\tn{val } V$ is the valence of the vertex $V$ of any curve $ \Gamma $
  of the given type. The one-dimensional cones of $\M_{0,n}$ are generated by
  the vectors $ v(I) $ of Construction \ref {moduli} for all $ I \subset
  \{1,\dots,n\} $ with $ 2 \le |I| \le n-2 $; by construction these are just the
  distance vectors of curves $ (\Gamma;x_1,\ldots,x_n)$ in $\M_{0,n}$
  having exactly one bounded edge of length $1$, with the markings $I$ to one
  and $ I^c $ to the other side. Hence there is only a finite number of
  combinatorial types of curves in $\M_{0,n}$. Note that the fan structure of
  $\M_{0,n}$ as in Construction \ref{moduli} coincides with the fan structure
  of the quotient of $\mathcal{G}(2,n)$ by its lineality space, see
  \cite[Theorem 4.2 and proof of Theorem 3.4]{SS}.
\end {remark}

\begin{remdef} [Notions of intersection theory for our purpose, see \cite {AR}
    for more details] \label {intersection}
  A \textit{tropical (affine) $k$-cycle} is a weighted fan of pure
  dimension $k$ in $V$ satisfying the balancing condition but where the range of
  the weight function is $\Z$. Hence the difference between a tropical fan of
  dimension $k$ and a tropical $k$-cycle is just the range of the weight
  function. A tropical $k$-cycle in a tropical fan $X$ is simply a $k$-cycle
  that is a subfan of $X$. Typically a tropical $k$-cycle in $X$ is an
  equivalence class of weighted fans of pure dimension $k$, where two fans are
  equivalent if they have a common refinement \cite[Definition 2.8]{AR}.
  However, in this paper we will not distinguish between equivalence
  classes and representatives of the equivalence class as the difference does
  not concern us here.

  As before, we will denote by $ \bet{X} $ the subset of $V$ of all cones of $X$
  (with non-zero weight). A \textit{(Weil) divisor} $D$ on $X$ is a
  cycle in $X$ of codimension 1.

  Let $X$ be a $k$-cycle. A (non-zero) \textit{rational function on} $X$ is a
  continuous function $ f: \bet{X} \to \R $ which is integer affine linear on
  each cone $ \sigma \subset \bet{X} $. We denote the linear part of the
  restriction of $f$ to $ \sigma $ by $ f_\sigma $.

  Let $X$ be a $k$-cycle and $f: \bet{X} \to \R$ a rational function on $X$.
  The \textit{Weil divisor $D(f)$ associated to} $f$ is the divisor in $X$
  consisting of the codimension-$1$ cones $ \tau \in X^{(k-1)} $ with weights
    \[ w(\tau) = \sum_{\sigma>\tau} w(\sigma) \, f_\sigma (v_{\sigma/\tau})
               - f_\tau \Big(
                 \sum_{\sigma>\tau} w(\sigma) \, v_{\sigma/\tau}
                 \Big), \]
  where the sum runs over all cones $\sigma$ s.t.\ $\tau \subsetneq \sigma$ is a face of
  $\sigma$, and $v_{\sigma/\tau} $ denotes the primitive normal vector of $ \sigma $
  relative to $ \tau $. It has been shown in \cite[Proposition 3.7]{AR} that this is
  indeed a tropical cycle.
  \end{remdef}

\begin{remdef} [Tropical Psi-classes $\psi_i$] \label {psi}
  Fix $n >2$ and $i \in \lk 1,\ldots,n\rk$. The \textit{Psi-class} $\psi_i
  \subset \M_{0,n}$ is the subfan of $ \M_{0,n} $ consisting of all cones of
  $\M_{0,n}$ of curves $(\Gamma; x_1,\ldots,x_n)$ such that the marked edge $ x_i$ is
  adjacent to a vertex of
  valence at least $4$ \cite{M}. Giving each top-dimensional cone the weight
  $1$ it has the structure of a tropical subfan of $\M_{0,n}$ of codimension 1.
  Using the language of \cite{AR}, we can rephrase this as: $\psi_i$ is a
  tropical Weil divisor associated to a rational function as proven in
  \cite[Proposition 3.5]{KM}. Note that a Psi-class is not defined up to
  rational equivalence as it is the case in classical geometry (and should
  therefore better be named Psi-divisor).  
\end{remdef}

\begin {convention} \label{conv-sk}
  When we talk in the following of ``the'' codimension-$k$ skeleton of $
  \M_{0,n}$ respectively $\psi_i$, we mean this w.r.t.\ the fan structure of
  Remarks \ref{combtypes} and \ref{psi}. 
\end {convention}

\begin{remark} [Vital divisors]
  Going back to Sean Keel \cite{K}, in the classical world the vital divisors
  generate the divisor class group of the moduli spaces $\overline{M}_{0,n}$ of
  rational $n$-marked stable curves and are sometimes called \textit{boundary
  divisors of} $\overline{M}_{0,n}$. The tropical analogue has first been
  studied by Matthias Herold in his diploma thesis \cite{H}. The results are
  also partly mentioned in \cite[Definition 2.4 and Lemma 2.5]{R}.
  We will introduce them here carefully.
\end{remark}

\begin{definition} [Tropical vital divisors in $\M_{0,n}$] \label {vital}
  Let $ n \ge 4 $ and $ S \subset \lk 1,\ldots,n \rk $ with $ 2 \le \bet{S} \le
  n-2 $. As $ \M_{0,n} $ is a simplicial fan with rays generated by the vectors
  $ v(I) $ of Remark \ref {combtypes} there is a unique piecewise linear
  function $ \phi_S $ on $ \M_{0,n} $ linear on each cone that takes the value
  $1$ on the vector $ v(S) $ and $0$ on all other rays. The associated Weil
  divisor $ D(\phi_S) $ to this function is called \textit{vital divisor
  w.r.t.\ $S$} and is denoted by $ D^S $.
\end{definition}

As $ v(S)=v(S^c) $ it is clear from the definition that $ D^S = D^{S^c} $.
Moreover, the following proposition taken from \cite[proof of Theorem 3.7]{H}
computes the weight function of the divisor $ D^S $ explicitly:

\begin{proposition} [Weight functions of the vital divisors] \label {vital-weight}
  Let $ n \ge 4 $ and $ S \subset \{1,\dots,n\} $ with $ 2 \le |S| \le n-2 $.
  Let $ \tau \in \M_{0,n}^{(n-4)} $ be a cone of codimension $1$; by
  Remark \ref {combtypes} its corresponding curves have exactly one
  $4$-valent vertex. Denote by $ A,B,C,D \subset \{1,\dots,n\} $ the sets
  of marked ends behind each of the four edges adjacent to this vertex. Then
  the weight of the vital divisor $ D^S $ on the cone $ \tau $ is
    $$ \qquad\quad
       w(\tau) = \left\{ \begin{array}{ll}
         1 & \tn{ if } S=I \cup J, \tn{ where }
           I \neq J \tn{ and } I,J \in \lk A,B,C,D \rk, \\
         -1 & \tn{ if } S=I \tn{ or } S^c=I \tn{ with }
           I \in \lk A,B,C,D \rk, \\
         0 & \tn{ otherwise}.
       \end{array} \right. $$   
\end{proposition}

\begin{proof}
  Adjacent to our given codimension-$1$ cone $ \tau $ of $ \M_{0,n} $ there are
  three top-dimensional cones $ \sigma_1,\sigma_2,\sigma_3 $ corresponding to
  resolving the $4$-valent vertex in the three ways shown in the picture below.

  \begin {center} \input {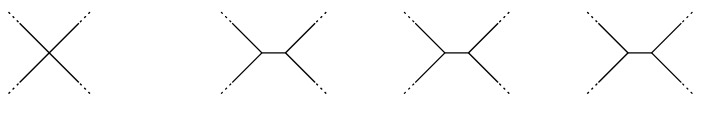} \end {center}
	
  The primitive normal vectors $ v_{\sigma_i/\tau} $ for $ i=1,2,3 $ as in
  Remark \ref {intersection} are then just the distance vectors of the three
  curves shown above on the right, with the newly added edge having length $1$
  and all other edges shrunken to zero, i.e.\ we have
    \[ v_{\sigma_1/\tau} = v(A \cup B), \quad
       v_{\sigma_2/\tau} = v(A \cup C), \mbox { and} \quad
       v_{\sigma_3/\tau} = v(A \cup D) \]
  with the vectors $ v(\,\cdot\,) $ as defined in Remark \ref {combtypes}. As
  $\M_{0,n}$ is balanced \cite[proof of Theorem 3.7]{GKM} we know moreover that
    \[ v(A \cup B) + v(A \cup C) + v(A \cup D) = v(A) + v(B) + v(C) + v(D). \]
  Using Definition \ref {vital} of the function $ \phi_S $, an application of
  the formula of Definition \ref {intersection} for the weight $ w(\tau) $ of
  the divisor $ D^S = D(\phi_S) $ on $ \tau $ thus yields (note that the
  weights of $ \M_{0,n} $ are all $1$)
  \begin{align*}
    w(\tau) &= \sum_{i=1}^3 (\phi_S)_{\sigma_i} (v_{\sigma_i/\tau})
               - (\phi_S)_\tau \Big(
                 \sum_{i=1}^3 v_{\sigma_i/\tau}
               \Big) \\
            &= \phi_S(v(A \cup B))+\phi_S(v(A \cup C))+\phi_S(v(A \cup D))
             - \phi_S(v(A)+v(B)+v(C)+v(D)) \\
            &= \begin{cases}
               1 & \tn{ if } S=I \cup J \tn{ with } I \neq J \mbox { and } I,J
                   \in \lk A,B,C,D \rk, \\
               -1 & \tn{ if } S=I \tn{ or } S^c=I \tn{ with } I \in \lk
                   A,B,C,D \rk, \\
               0 & \tn{otherwise}.
               \end{cases}
  \end{align*}
\end{proof}

\begin{example} [Vital divisors in $\M_{0,5}$] \label {ex-vital-5}
  By symmetry, in $ \M_{0,5} $ we just have to consider sets $S$ with
  $\bet{S}=2$. There are $10$ such sets and thus $10$ vital divisors. Set for example
  $S=\lk 1,2 \rk$. Curves of codimension 1 in $ \M_{0,5} $ have one bounded
  edge on exactly one of whose ends lies a 4-valent vertex. So in the notation
  of Proposition \ref {vital-weight} we have w.l.o.g.\ $\bet{A}=2$;
  $A$ determines the type of the curve $\M_{0,5}$ uniquely. Using the
  weight function above, we get the following weights:
  \begin{center}
    \begin{tabular}{|c|c@{\;\;}c@{\;\;}c@{\;\;}c@{\;\;}c@{\;\;}c@{\;\;}%
                    c@{\;\;}c@{\;\;}c@{\;\;}c|} \hline
      $A=$ & $\lk 1,2\rk$ & $\lk 1,3 \rk$ & $\lk 1,4 \rk$ & $\lk 1,5 \rk$ &
             $\lk 2,3\rk$ & $\lk 2,4 \rk$ & $\lk 2,5 \rk$ & $\lk 3,4 \rk$ &
             $\lk 3,5\rk$ & $\lk 4,5 \rk$ \\ \hline
      weight & $-1$&$0$&$0$&$0$&$0$&$0$&$0$&$1$&$1$&$1$ \\ \hline
    \end{tabular}
  \end{center} \hfill $\circ$
\end{example}

  We now want to check that codimension-$k$ skeletons of and Psi-classes in
$\M_{0,n}$ are in fact tropical fans, i.e.\ that they satisfy the balancing
condition. So let us fix $ k \in \lk 0,\ldots, n-3\rk $ and a weight function
$w: \M_{0,n}^{(n-3-k)} \to \Z_{>0}$ on the codimension-$k$ skeleton of $
\M_{0,n} $. Moreover, let $ \tau \in \M_{0,n}^{(n-3-k-1)}$ be a cone,
corresponding by Remark \ref {combtypes} to a certain combinatorial type of
$n$-marked curves. In order to verify the balancing condition for $w$ at $ \tau
$, the following toolkit will be useful.

\begin{remark} [About the balancing condition in $ \M_{0,n} $]
    \label{splitting}
  The cones $\sigma \in \M_{0,n}^{(n-3-k)}$ containing $\tau$ can be obtained
  by resolving one chosen vertex $V$ of valence $r$ at least $4$ in $ \tau $ in
  the same way as it can be resolved in $ \M_{0,r} $ by adding one
  bounded edge. More precisely, this means that we replace $V$ by two vertices
  joined by a bounded edge, with the $r$ edges of $V$ split up in every
  possible way onto the two new vertices such that there are at least two such edges
  on each side of the bounded edge. The following picture shows these types of
  resolutions; it can either be seen as a picture of curves in $ \M_{0,r} $ or
  as a \textit {local} picture of curves in $ \M_{0,n} $ around $V$.

  \begin {center} \input {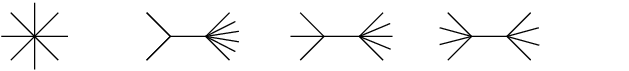} \end {center}

  In order to check the balancing condition at $ \tau $ it therefore suffices
  to split the total sum $\sum_{\sigma > \tau} w(\sigma) \, v_{\sigma/\tau}$
  into parts, where each part corresponds to the resolution of one vertex $V$,
  and verify it for each part separately. This leads to the following lemma
  which states that the balancing condition for codimension-$k$ cells in $
  \M_{0,n} $ can be split up into several conditions in lower-dimensional
  moduli spaces.
\end{remark}

\begin{lemma} [Splitting the balancing condition] \label{local}
  With notations as above, let $\lk V_1,\ldots,V_m \rk$ be the set of $4$- or
  higher valent vertices of a curve in $\tau$, and let $ r_1,\dots,r_m
  $ be their respective valences. Then to verify the balancing condition
  at $ \tau $ in $ \M_{0,n} $ it suffices to check the balancing at all
  $0$-dimensional cones $ \tau_i $ in $ \M_{0,r_i}$ for $i \in \lk 1,\ldots,
  m\rk $ (corresponding to curves (``stars'') having only one vertex of valence
  $r_i$).

  \begin {center} \input {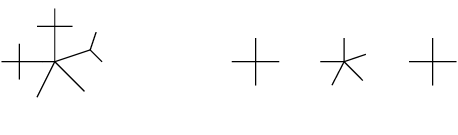} \end {center}
\end{lemma}

\begin{proof}
  Consider a 4- or higher valent vertex $ V \in \{ V_1,\dots,V_m \} $ of $ \tau
  $, let $r$ be its valence, and denote by $ A_i $ for $ i=1,\dots,r $ the set
  of marked edges behind the $i$-th edge of $V$. Thus we have $ \sqcup_i A_i= \lk
1,\ldots,n \rk$. Define a linear map $\phi$ by
    \[ \phi: \R^{\binom r2} \to \R^{\binom n2}, \quad
         (x_{i,j})_{i<j} \mapsto (\tilde x_{k,l})_{k<l} \]
  where
    \[ \tilde x_{k,l} = \begin {cases}
                          x_{i,j} &
                            \mbox {if $ (k,l) \in A_i \times A_j $ or
                                      $ (l,k) \in A_i \times A_j $ for some
                                      $ 1 \le i < j \le r $,} \\
                          0 & \mbox {otherwise} \\
                        \end {cases} \]
  for $ 1 \le k < l \le n $. Then by construction the distance vectors of
  Construction \ref {moduli} are transformed by $ \phi $ as
    \[ \phi ( \tilde v(I)) = \tilde v \Big( \bigcup_{i \in I} A_i \Big)
       \tag {1} \]
  for all subsets $ I \subset \{ 1,\dots,r \} $. In particular, for $ i=1,
  \dots,r $ the $i$-th basis vector $ d_i $ of the lineality space of
  $\M_{0,r}$ is mapped to $ \phi(d_i)=\tilde{v}(A_i) $.

  Let $ \tau_\star $ be the $0$-dimensional cone in $ \M_{0,r} $,
  corresponding to the star curve with one vertex. By Remark \ref {splitting}
  a $1$-dimensional cone $ \sigma_\star $ around $ \tau_\star $ corresponds to
  a cone $ \sigma $ around $ \tau $ in $ \M_{0,n} $ describing the same local
  resolution. Let us assume that the balancing condition holds at $
  \tau_\star $, i.e.\ that
    \[ \sum_{\sigma_\star > \tau_\star}
         w(\sigma) \, \tilde v_{\sigma_\star/ \tau_\star}=\sum_i a_i \, d_i
       \qquad \in \R^{\binom r2}
       \tag {2} \]
  for some $ a_i \in \R $. By Remark \ref {combtypes} the normal vectors
  $ \tilde v_{\sigma_\star/ \tau_\star} $ in $ \M_{0,r} $ are exactly $ \tilde
  v(I) $ for the corresponding subset $ I \subset \{ 1,\dots,r \} $ with $ 2
  \le |I| \le r-2 $, and by (1) these vectors are mapped by $ \phi $ to the
  corresponding normal vectors $ \tilde v_{\sigma/\tau} $ in $ \M_{0,n} $. So
  applying $ \phi $ to (2) we get
    \[ \sum_{\sigma_* > \tau_*}
         w(\sigma) \, \tilde v_{\sigma/\tau}=\sum_i a_i \, \tilde v(A_i)
       \qquad \in \R^{\binom n2}. \]
  The vectors $ \tilde v(A_i) $ lie in the lineality space for $ |A_i|=1 $
  and in $ V_\tau $ otherwise, so taking the quotient by these spaces this sum
  reduces to zero. The claim of the lemma thus follows with the second part of
  Remark \ref {splitting}.
\end{proof}

\begin{proposition} [Codimension-$k$ skeleton of $\M_{0,n}$] \label {comod}
  Let $k \in \lk 0,\ldots, n-3 \rk$. Then the codimension-$k$ skeleton
  ${\M_{0,n}^{(n-3-k)}}_\tn{sk}$ of $\M_{0,n}$ with the weight function $
  w:\M_{0,n}^{(n-3-k)} \to \lk 1\rk $ is balanced. Hence, it is a tropical fan.
\end{proposition}

\begin{proof}
  By Lemma \ref {local} we can reduce the proof to the local situation of a
  star curve with a vertex of some valence $ r \in \{ 4,\dots,k+3 \} $, i.e.\
  to the balancing condition around the vertex $ \tau_\star $ in $ \M_{0,r} $.

  Let us think of this balancing condition in terms of coordinate vectors
  in $ \R^{\binom r2} $. We have to compute the sum $ \tilde v $ of all normal
  vectors $ \tilde{v}_{\sigma_\star / \tau_\star} $ arising from resolving the
  vertex of the star. Consider the first entry of this vector: here, we sum up
  $1$ a number of times, a $1$ for each type where the marked edges $1$ and $2$
  lie on opposite sides of the bounded edge. By symmetry (i.e.\ no marked edge
  is distinguished), this sum is the same in each other entry of $ \tilde v $.
  So $ \tilde v $ is a multiple of the vector $ (1,\ldots,1)^\top $. Consider
  now the lineality space: summing up all vectors $ d_1,\dots,d_r $ gives $ 2
  \cdot (1,\ldots,1)^\top$. So $ \tilde v $ is the zero vector modulo the
  lineality space.
\end{proof}

\begin{proposition} [Codimension-$k$ skeleton of a Psi-class $\psi_i$ in
    $\M_{0,n}$] \label {copsi}
  Let $k \in \lk 0,\ldots, n-4 \rk$ and $i \in \lk 1,\ldots,n \rk$. Then
  the codimension-$k$ skeleton ${\psi_i^{(n-4-k)}}_{\tn{sk}}$ of the $i$-th
  Psi-class $\psi_i$ of $\M_{0,n}$ is balanced for the weight function
  $w:\psi_i^{(n-4-k)} \to \lk 1\rk$. Hence, it is a tropical fan.   	
\end{proposition}

\begin{proof}
  Again we can use Lemma \ref {local} to reduce the proof to the local
  situation of a star with a vertex of some valence $ r \in \{ 4,\dots,k+4 \}
  $. If the $i$-th unbounded edge is not adjacent to the chosen
  vertex then the computation is exactly the same as in the proof of
  Proposition \ref {comod}, so let us assume that it is adjacent to the chosen
  vertex. Then, as in the picture of Example \ref {ex-psi} below, we only have
  to consider resolutions of the star in which the $i$-th edge remains adjacent
  to a vertex of valence at least $4$.

  Think again in terms of coordinate vectors in $ \M_{0,r} $. Again, denote the
  sum of all normal vectors of these resolutions by $ \tilde v $. Then the $
  (j,k) $-coordinates of $ \tilde v $ with $ j,k \neq i $ are all the same by
  symmetry, and likewise for the $ (j,k) $-coordinates where $ j=i $ or $ k=i
  $. So if we set w.l.o.g.\ $ i=1 $ then we can write $ \tilde v $ as
    \[ \tilde v = (\underbrace{M,\ldots,M}_{r-1 \tn{ times }},N,\ldots,N)^\top
                = M \, d_1 + \frac N2 \, (-d_1+d_2+\cdots+d_r) \]
  for suitable $ M,N \in \N $. As this is a vector in the lineality space, the
  balancing condition follows.
\end{proof}

\begin {example} [Codimension-$2$ skeleton of $\psi_1$ in $\M_{0,7}$]
    \label {ex-psi}
  The codimension-2 skeleton of $ \psi_1 $ in $ \M_{0,7} $ is 1-dimensional,
  and thus there is only one balancing condition to check, namely that around
  the 0-dimensional cell corresponding to the star curve with only one vertex.
  Allowed resolutions of this type are:

  \begin {center} \input {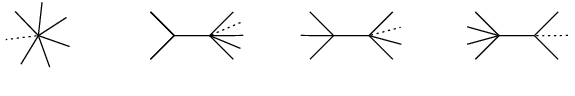} \end {center}

  As in the proof of Proposition \ref {copsi} let $ \tilde v $ be the sum of
  all normal vectors of these resolutions. Write $ \tilde v $ as $ \tilde v =
  \tilde v_{\textrm {(a)}} + \tilde v_{\textrm {(b)}} + \tilde v_{\textrm
  {(c)}} $, where the summands denote the parts of $ \tilde v $ arising from
  resolutions of type (a), (b), (c) as in the picture above, respectively. Note
  that each such type corresponds to various resolutions corresponding to the
  choice of labeling of the marked ends. In the array below, for each
  resolution type the coordinates of $ \tilde v $ are listed, where $ j,k \neq
  1 $. 
  \begin{center}
    \begin{tabular}{|c||c|c|c|c|} \hline
      coordinate & $ \tilde v_{\textrm {(a)}} $ &
                   $ \tilde v_{\textrm {(b)}} $ &
                   $ \tilde v_{\textrm {(c)}} $ & $ \tilde v $ \\ \hline \hline
      $ (1,j) $ & $\binom{5}{4}=5$ & $\binom{5}{3}=10$ & $\binom{5}{2}=10$ &
        $25$ \\	\hline
      $ (j,k) $ & $2 \cdot \binom{4}{3}=8$&$2 \cdot \binom{4}{2}=12$&$2
        \cdot \binom{4}{1}=8$ &$28$ \\ \hline
    \end{tabular}
  \end{center}
  For example, the $(1,j)$-coordinate $10$ of $ \tilde v_{\textrm {(b)}} $
  corresponds to the 10 choices of distributing the remaining labels on the
  ends in type (b) if the unbounded edge $j$ has been put at the left vertex.

  So we have
    \[ \tilde v_{\textrm {(a)}} = 4 \sum_{j=1}^7 d_j - 3 d_1, \qquad
       \tilde v_{\textrm {(b)}} = 6 \sum_{j=1}^7 d_j - 2 d_1, \qquad
       \tilde v_{\textrm {(c)}} = 4 \sum_{j=1}^7 d_j + 2 d_1. \]
  Note that it is not just the sum $ \tilde v $ that is zero modulo the
  lineality space, but also the individual vectors $ \tilde v_{\textrm {(a)}}
  $, $ \tilde v_{\textrm {(b)}} $, $ \tilde v_{\textrm {(c)}} $ corresponding
  to the resolution types themselves. In fact, the proof of Proposition \ref
  {copsi} shows that the analogous statement holds for the (one-dimensional)
  codimension-$ (n-5) $ skeleton of a Psi-class in $ \M_{0,n} $ for all $ n \ge
  5 $ since the symmetry argument given there also applies if we only consider
  a single resolution type. \hfill $\circ$
\end {example}

\begin {remark} [Choice of weight function]
  Of course, the computations above depend on the chosen weight function. The
  moduli spaces $ \M_{0,n} $ and the Psi-classes $ \psi_i $, considered as
  tropical cycles, have all weights of their facets equal to $1$, and thus
  it was natural in Propositions \ref {comod} and \ref {copsi} to also equip
  the $k$-skeletons of these cycles with the constant weight function 1. For
  other cycles such as for instance intersection products $ \psi_1^{m_1} \cdot \ldots
  \cdot \psi_n^{m_n} $ of Psi-classes (where $ m_1,\ldots,m_n \in \N_0 $)
  this is in general no longer the case, especially when they are of codimension $2$ or
  higher. 
\end {remark}

  Now we want to explore the question of uniqueness of the weight function for
certain tropical subfans of $\M_{0,n}$. This question is strongly related to
the one of irreducibility introduced in \cite {GKM}. To simplify notations, we
will assume from now on that all weights of the cones in tropical cycles are
non-zero.

\begin {definition} [Irreducible tropical cycles] \label {def-irred}
  A tropical cycle $X$ in $ V=\Lambda \otimes_\Z \R $ is said to be
  \textit{(globally) irreducible} if there does not exist a tropical cycle $Y$
  of the same dimension in $V$ such that $ |Y| \subsetneq |X| $. Of course,
  this definition then applies to tropical fans as well.
\end {definition}

\begin{remark} \label{rem-class}
  Lemma \ref{criterion} implies that, just as in the classical situation, the
  support of a tropical cycle can always be written as the union of the
  supports of irreducible tropical cycles. However, such a decomposition is in
  general not unique \cite[Remark 2.19]{GKM}.
\end{remark}

\begin {remark} \label {rem-irred}
  In Definition \ref {def-irred} the cones of $Y$ are just required to be
  contained in the union $ |X| $ of the cones of $X$. However, the definition
  does not change if one requires all cones of $Y$ to actually be cones of $X$,
  just with possibly different weights. To see this, assume there is a tropical
  cycle $Y$ satisfying $ |Y| \subsetneq |X|$. By passing to a common refinement
  with $X$, we can then first of all make sure that every cone of $Y$ is
  contained in a cone of $X$. But then all cones of $Y$ contained in the same
  cone of $X$ must have the same weight due to the balancing condition, and
  hence can be made into a single cone.
\end {remark}

\begin{lemma} \label {criterion}
  A tropical cycle $X$ is irreducible if and only if ``its weight function
  is unique up to a global multiple'', i.e.\ if and only if for every cycle $Y$
  of the same dimension and consisting of at most the cones of $X$ there is a
  rational number $ \lambda \in \Q $ such that $ w_Y(\sigma) = \lambda \,
  w_X(\sigma) $ for every cone $ \sigma $ of $X$.
\end{lemma}

\begin{proof}
  ``$ \Rightarrow $'': This is \cite[Lemma 2.21]{GKM}.

  ``$ \Leftarrow $'': Let $Y$ be a cycle with $ |Y| \subsetneq |X| $. By Remark
  \ref {rem-irred} we can assume that each cone of $Y$ is a cone of $X$, so
  there must be a cone $ \sigma $ of $X$ with $ w_Y(\sigma) = 0 $. But this
  requires $ \lambda=0 $ in our assumption, so $Y$ would have to be the zero
  cycle.
\end{proof}

Using this lemma the problem of checking irreducibility can often be reduced to
a local computation. For this we need the concept of a cycle that is connected
in codimension 1.

\begin {definition} [Locally irreducible cycles, connected in codimension 1]
    \label {def-locirr}
  Let $X$ be a tropical $k$-cycle.
  \begin {enumerate}
  \item \label {def-locirr-a}
    The cycle $X$ is called \textit {locally irreducible} if for each $
    (k-1) $-dimensional cone $ \tau $ there is up to a common multiple only
    one choice of weights $ w(\sigma) $ for the cones $ \sigma $ adjacent to $
    \tau $ such that the balancing condition $ \sum_{\sigma>\tau} w(\sigma) \,
    v_{\sigma/\tau} = 0 \mod V_\tau $ holds at $ \tau $, see \cite[Section 1.6]{R}.
  \item \label {def-locirr-b}
    The cycle $X$ is called \textit{connected in codimension 1} if any two
    $k$-dimensional cones $ \sigma,\sigma' $ in $X$ can be connected by a
    sequence $ \sigma=\sigma_0,\ldots,\sigma_m=\sigma'$ of $k$-dimensional
    cones in $X$ such that $ \sigma_{i-1} \cap \sigma_i $ is of dimension $ k-1
    $ for all $ i=1,\dots,m $ \cite[Section $3$]{BJSST}.
  \end {enumerate}
\end {definition}

\begin{remark} \label{rem-crit}
  A criterion similar to Definition \ref{def-irred} holds also for locally
  irreducible cycles: we just have to apply Lemma \ref{criterion} to the stars
  corresponding to each $\tau$ to see that a $k$-dimensional cycle is locally
  irreducible if and only if its stars at all $ (k-1) $-dimensional cones are
  irreducible in the sense of Definition \ref{def-irred}.
\end{remark}

\begin{lemma} [Criterion for irreducibility of a tropical cycle]
    \label {uniqueness}
  Every locally irreducible tropical cycle that is connected in codimension 1
  is also globally irreducible.
\end{lemma}

\begin{proof}
  Let $ \sigma, \sigma'\in X^{(k)} $ be two top-dimensional cones in a
  $k$-dimensional cycle $X$. As $X$ is connected in codimension 1 there exist
  cones $ \sigma_i \in X^{(k)}$ for $i \in \lk 1,\ldots,m \rk$ forming a
  sequence $ \sigma=\sigma_0,\ldots,\sigma_m=\sigma'$ such that each
  intersection $ \sigma_{i-1} \cap \sigma_i $ is a cone of dimension
  $k-1$. As $X$ is locally irreducible, the weight $ w(\sigma_i) $ determines
  the weight $ w(\sigma_{i-1}) $ uniquely. Hence, running through the sequence,
  we can conclude that $ w(\sigma) $ determines $ w(\sigma') $ uniquely and
  thus that $X$ is globally irreducible by Lemma \ref {criterion}.
\end{proof}

The converse of this lemma does not hold, see Example \ref {ex-m06}.

\begin{lemma} \label{mcon}
  $\M_{0,n}$ is connected in codimension $1$ for each $n >2$.
\end{lemma}

\begin{proof}
  Let $ \sigma $ be a top-dimensional cone of $ \M_{0,n} $, corresponding by
  Remark \ref {combtypes} to a combinatorial type of 3-valent tropical curves.

  First we show that $ \sigma $ can be connected in codimension 1 (as in
  Definition \ref {def-locirr} \ref {def-locirr-b}) to a top-dimensional
  cone corresponding to a \textit {caterpillar tree}, i.e.\ to a type of 3-valent curves
  for which the number $m$ of 3-valent vertices that are not adjacent to any
  end is zero. So assume that $ m>0 $. As in the picture below we can then pick
  such a vertex $V$ that is not adjacent to any end but adjacent to at least
  one further vertex $W$ with two ends attached to it.

  \begin {center} \input {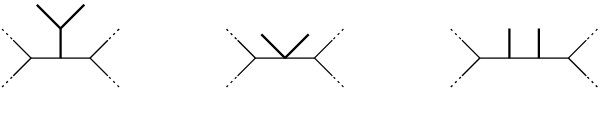} \end {center}
  
  The picture then shows how $ \sigma $ is connected in codimension 1 to
  another type $ \sigma' $ for which the value of $m$ is reduced by 1 --- by
  contracting the edge joining $V$ and $W$ to a 4-valent vertex, and then
  resolving this 4-valent vertex again in a different way. Continuing this
  process, we can connect $ \sigma $ in codimension 1 to a caterpillar tree.

  But now any two caterpillar trees are connected in codimension 1 by a sequence of
  exchanges of neighboring ends as in the following picture:

  \begin {center} \input {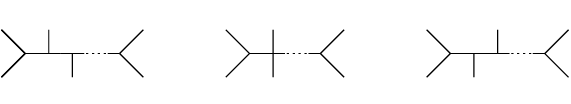} \end {center}

  Hence we conclude that any two top-dimensional cones of $ \M_{0,n} $ are
  connected in codimension 1.
\end{proof}

\begin{remark}
 The statement of Lemma \ref{mcon} is not new. It is noticed in
 \cite[Section $3$]{SS} that the Pl\"{u}cker ideal $I_{2,n}$ is prime and that
 $\mathcal{G}(2,n)$ is the tropical variety defined by $I_{2,n}$. By
 \cite[Theorem 3.1]{BJSST} it follows that  $\mathcal{G}(2,n)$ is connected in
 codimension one. When taking the quotient by its lineality space in order to
 obtain $\M_{0,n}$ this property is conserved.   
\end{remark}

We are now going to prove that for all $n>2$ the tropical moduli space
$\M_{0,n}$ and its vital divisors, respectively Psi-classes are irreducible.
These results also hold in classical algebraic geometry except for Psi-classes
which are not divisors but only divisor classes: for $\overline{M}_{0,n}$ one
can use the quotient description of \cite[Theorem 6.1]{K2}, and the classical
vital divisor $D^S$ is isomorphic to the product $\overline{M}_{0,|S|+1}\times
\overline{M}_{0,|S^c|+1}$ of irreducible varieties \cite[Corollary 3.9]{K1},
and hence is irreducible itself.

\begin{proposition} \label {mir}
  $\M_{0,n}$ is irreducible for each $n>2$.
\end{proposition}

\begin{proof}
  By Lemmata \ref{uniqueness} and \ref{mcon} it suffices to check local
  irreducibility. We have to look at the weights of cones corresponding to the
  resolving curves of a curve with one $4$-valent vertex. By Lemma \ref {local}
  we can assume that this curve is just a star with one $4$-valent vertex in
  the middle, living in $\M_{0,4}$.

  \begin {center} \input {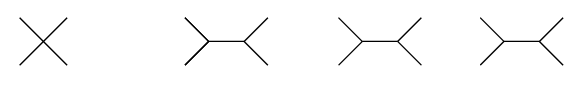} \end {center}

  So the task translates to finding non-trivial solutions $(a_1,\ldots,a_7)$
  for the linear system of equations
    $$  a_1 \, \tilde{v}(\lk 1,2\rk)
      + a_2 \, \tilde{v}(\lk 1,3\rk)
      + a_3 \, \tilde{v}(\lk 1,4\rk)
      = a_4 \, d_1 + a_5 \, d_2 + a_6 \, d_3 + a_7 \, d_4 $$
  in $ \R^{\binom 42} = \R^6 $. Solving this system gives $a_i=1$ for all $i\in
  \lk 1,\ldots,7\rk$ up to a multiple. Hence $ \M_{0,n} $ is locally
  irreducible and thus also globally irreducible.
\end{proof}

Next, we want to study the irreducibility of codimension-1 cycles in $ \M_{0,n}
$ --- in particular of Psi-classes and vital divisors. In this case we can
often not use Lemma \ref {uniqueness} as the cycles in question are not locally
irreducible. Let us start by illustrating this fact with an example that will
also be needed later.

\begin {example}[Difference between local and global irreducibility]
    \label {ex-m06}
  The following picture shows a type of curves corresponding to a
  codimension-2 cell $ \tau $ in $ \M_{0,6} $ with the types of its six
  adjacent codimension-1 cells $ \sigma_1,\dots,\sigma_6 $:

  \begin {center} \input {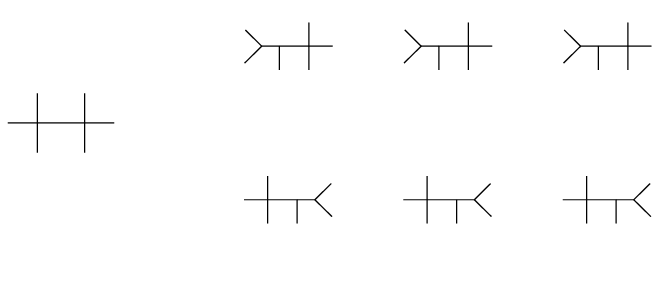} \end {center}

  By Lemma \ref {local} and the balancing condition in $ \M_{0,4} $ we have
  $ v_{\sigma_1/\tau} + v_{\sigma_2/\tau} + v_{\sigma_3/\tau} = 0 $ and
  $ v_{\sigma_4/\tau} + v_{\sigma_5/\tau} + v_{\sigma_6/\tau} = 0 $. In fact,
  an easy computation in coordinates shows that the linear equation
  $ \sum_{i=1}^6 a_i \, v_{\sigma_i/\tau} = 0 $ is satisfied in $ \M_{0,6} $ if
  and only if $ a_1=a_2=a_3 $ and $ a_4=a_5=a_6 $. All six resolutions above
  contribute to the vital divisor $ D^{\{1,2,3\}} $ with the same non-zero weight $-1$ by
  Proposition \ref {vital-weight}. This means that we can scale $ a_1,a_2,a_3 $
  by one multiple and $ a_4, a_5, a_6 $ independently by another one in order to
  satisfy the equation $ \sum_{i=1}^6 a_i \, v_{\sigma_i/\tau} = 0 $. Therefore $
  D^{\{1,2,3\}} $ is not locally irreducible
  at $ \tau $. We will see in Proposition \ref {prop-vital} however that vital
  divisors are always globally irreducible. \hfill $\circ$
\end {example}

\begin {lemma} \label {lem-divisor-set}
  Let $Y$ be a codimension-1 cycle in $ \M_{0,n} $. Let $ \sigma $ be a
  top-dimensional cone of $Y$; by Remark \ref {combtypes} its corresponding
  curves have exactly one $4$-valent vertex. Denote by $ A,B,C,D \subset
  \{1,\dots,n\} $ the sets of marked ends behind each of the four edges
  adjacent to this vertex (in any order). Then the weight $ w(\sigma) $ in $Y$
  depends only on the sets $ A,B,C,D $ (and not on the subtrees behind the
  $4$-valent vertex).
\end {lemma}

\begin {proof}
  Let $ \tau' $ be a codimension-1 cell of $Y$, i.e.\ a codimension-2 cell in $
  \M_{0,n} $. By Remark \ref {combtypes} the corresponding curves have either
  two 4-valent vertices or one 5-valent vertex, with all other vertices being
  3-valent. The case of a 5-valent vertex is irrelevant for the claim of the lemma and
  will be treated in Proposition \ref{pir}. Let us consider the case when $ \tau' $ has
two  4-valent vertices.
  Without loss of generality we may then assume that the ends are labeled in
  such a way that there is a forgetful map $ \M_{0,n} \to \M_{0,6} $ sending $
  \tau' $ to the cone $ \tau $ in Example \ref {ex-m06} above. By Remark \ref
  {splitting} the top-dimensional cones $ \sigma'_1,\dots,\sigma'_6 $ in $Y$
  adjacent to $ \tau' $ are in one-to-one correspondence with the cones $
  \sigma_1,\dots,\sigma_6 $ in $ \M_{0,6} $ of Example \ref {ex-m06}. Now let $
  a_i $ be the weight of $ \sigma'_i $ in $Y$ for $ i=1,\dots,6 $. By the
  balancing condition we have $ \sum_i a_i \, v_{\sigma'_i/\tau_i} = 0 $ in $
  \M_{0,n} $. Applying the (linear) forgetful map $ \M_{0,n} \to \M_{0,6} $ we
  conclude that $ \sum_i a_i \, v_{\sigma_i/\tau_i} = 0 $ in $ \M_{0,6} $ as
  well, which by Example \ref {ex-m06} means that $ a_1=a_2=a_3 $ and $
  a_4=a_5=a_6 $. In particular, the types $ \sigma'_1 $, $ \sigma'_2 $, $
  \sigma'_3 $ have the same weight in $Y$. Applying this to the labeling where
  the forgetful map sends the marked ends in $B$, $C$, $D$ to $ x_4 $, $ x_5 $,
  $ x_6 $, respectively, we see by the argument of the proof of Lemma \ref
  {mcon} that any two top-dimensional cells of $Y$ that just differ by the
  arrangement of the markings in $A$ behind one of the edges of the 4-valent
  vertex must have the same weight. As the same holds for $B$, $C$, $D$ as
  well, the lemma follows.
\end {proof}

\begin{lemma} \label {lem-psi-conn}
  $\psi_i$ is connected in codimension 1 for each $i \in \lk 1, \ldots,n \rk$.
\end{lemma}

\begin{proof}
  Let $ \sigma $ be a top-dimensional cone in $ \psi_i $. As in the picture
  below on the left it corresponds to a combinatorial type of tropical curves
  with exactly one 4-valent vertex $V$ adjacent to the marked point $ x_i $,
  all other vertices being 3-valent.

  By the proof of Lemma \ref {lem-divisor-set} we see that $ \sigma $ can be
  connected in codimension 1 to any other such type where the trees behind the
  edges of $V$ have been replaced by caterpillars, and where within each such
  caterpillar the marked points have been ordered in a given way. The picture
  below then shows that it is also possible to connect such a cone in codimension 1
  to the one where a marked point $ x_j $ with $ j \neq i $ changes its place
  from behind one of the edges of $V$ to another --- through a codimension-1
  cone in $ \psi_i $ corresponding to curves with a 5-valent vertex.

  \begin {center} \input {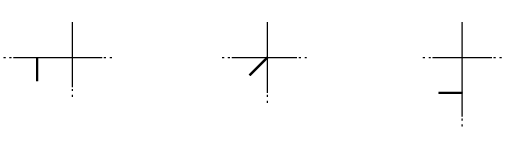} \end {center} 

  Using these steps we see that $ \sigma $ can be connected in codimension 1
  to any other top-dimensional cone of $ \psi_i $.
\end{proof}

\begin{proposition} \label{pir}
  Each Psi-class $\psi_i$ in $\M_{0,n}$ is irreducible. 
\end{proposition}

\begin{proof}
  Assume w.l.o.g.\ $i=1$. As in the proof of Proposition \ref{mir} it suffices
  by Lemmata \ref {uniqueness} and \ref {lem-psi-conn} to check local
  irreducibility at each codimension-1 cone of $ \psi_1 $. There
  are two types of these cones: those with two 4-valent vertices (one of them
  adjacent to $ x_1 $) and those with one 5-valent vertex (adjacent to $ x_1
  $). In the former case only the vertex that is not adjacent to $ x_1 $ can be
  resolved in $ \psi_1 $, and the computation to check local irreducibility is
  exactly the same as in Proposition \ref {mir}. In the latter case we again
  apply Lemma \ref {local} to reduce the situation to the resolutions of a
  star in $\M_{0,5}$ such that the marking $ x_1 $ remains adjacent to a
  $4$-valent vertex. There are six such resolving curves of the $5$-valent
  vertex, corresponding to a choice of two of the marked ends $ x_2,\dots,
  x_5 $ to put at the 4-valent vertex together with $ x_1 $ and the bounded
  edge. Hence, this time the linear system of equations reads
  \begin {align*}
    &   a_1 \, \tilde{v}(\lk 2,3 \rk)
      + a_2 \, \tilde{v}(\lk 2,4 \rk)
      + a_3 \, \tilde{v}(\lk 2,5 \rk)
      + a_4 \, \tilde{v}(\lk 3,4 \rk)
      + a_5 \, \tilde{v}(\lk 3,5 \rk)
      + a_6 \, \tilde{v}(\lk 4,5 \rk) \\
    & \qquad \qquad
      = \sum_{j=7}^{11} a_j \, d_{j-6}
  \end {align*}
  in $ \R^{\binom 52} = \R^{10} $. Solving this system gives up to a
  multiple $ a_j=1 $ for $ j=1,\dots,7 $ and $ a_j=2 $ for $ j=8,\dots,11 $ as
  the only solution. Hence $ \psi_i $ is locally irreducible and thus also
  globally irreducible.
\end{proof}

\begin {remark}
  Using similar arguments as in Lemmata \ref {mcon} and \ref {lem-psi-conn} one
  can show that the codimension-$k$ skeletons of the moduli spaces $ \M_{0,n} $
  as well as of the Psi-classes $ \psi_i $ (see Propositions \ref {comod} and
  \ref {copsi}) are connected in codimension 1 as well. However, they are in
  general not irreducible, as the following example shows.
\end {remark}

\begin{example} [Irreducibility of other tropical fans in $\M_{0,n}$]
  Consider the (one-dimen\-sional) codimension-$ (n-5) $ skeleton of the
  Psi-class $ \psi_1 $ in $ \M_{0,n} $ for $ n \ge 6 $, corresponding to
  resolutions of the star curve with one $n$-valent vertex that keep $ x_1 $ at
  a vertex of valence at least 4. We have seen at the end of Example \ref{ex-psi} 
 that in this case it is not only the sum of all normal vectors of
  the possible resolutions that sum up to zero modulo the lineality space, but
  also the vectors of each individual resolution type, i.e.\ when fixing how
  many marked ends should be adjacent to $ x_1 $. In particular, as there is
  more than one such resolution type if $ n \ge 6 $, this means that there are
  several choices of weights satisfying the balancing condition and thus that
  the codimension-$ (n-5) $ skeleton of $ \psi_1 $ in $ \M_{0,n} $ is not
  locally irreducible. As for this cycle there is only one balancing condition
  to check it follows that it is not globally irreducible either.

  A similar argument can be used to show that other codimension-$k$ skeletons
  of a Psi-class or of $ \M_{0,n} $ also fail to be locally irreducible.
  In general this does not imply however that these cycles are not globally
  irreducible either (e.g.\ we will see in the proof of Proposition \ref{prop-vital} 
  that the vital divisors of $ \M_{0,n} $ are globally but not
  locally irreducible). \hfill $\circ$
\end{example}

The following example will be of use in the remainder.

\begin{example} [Irreducible divisors in $\M_{0,5}$] \label {m05}
  In $\M_{0,5}$ the problem of finding irreducible divisors is equivalent
  to determine all minimal linear dependent sets of the 10 vectors $v(I)$
  where $I \subset \lk 1,\ldots,5\rk$ and $\bet{I}=2$, or in other words to
  find all circuits of the vector matroid given by these vectors. An explicit
  computation shows that there are two kinds of irreducible divisors: those
  containing $4$ cones in their support and those with 6 cones. Up to
  multiples, these divisors are the following:
  \begin {itemize}
  \item 4 cones: The only effective divisors with four cones are the cycles $
    \psi_i^\natural $ for $ i=1,\dots,5 $ given as the difference of the
    1-skeleton and $ \psi_i $ --- for example $\psi_1^\natural $ consists of the
    cones $v(\lk 1,2\rk)$, $v(\lk 1,3\rk)$, $v(\lk 1,4\rk)$, and
    $v(\lk 1,5\rk)$, each with weight 1. The only other divisors with four
    cones are the ten vital divisors $ D^S $ of $\M_{0,5}$. They have
    respectively three cones of weight $1$ and one of weight $-1$ (see
    Proposition \ref {vital-weight}).
  \item 6 cones: The Psi-classes $\psi_i$ are the only effective irreducible
    divisors with six cones. There are ten more irreducible divisors with six
    cones, given by $ \psi_i^\natural - \psi_j^\natural =: V^{\lk i, j\rk} $
    for $i>j$. \hfill $\circ$
  \end {itemize} 
\end{example}

\begin{proposition} \label {prop-vital}
  For $n>3$ all vital divisors $D^S$ in $\M_{0,n}$ are irreducible.
\end{proposition}

\begin{proof}
  The case $n=4$ is clear as the only vital divisor corresponds to the
  star with a $4$-valent vertex. So let $ n \ge 5 $, and let $ S \subset
  \{1,\dots,n\} $ be a subset with $ 2 \le |S| \le n-2 $, giving rise to a
  vital divisor $ D^S $ in $ \M_{0,n} $. By Proposition \ref {vital-weight}
  the cones occurring in $ D^S $ correspond to curves with one 4-valent vertex
  such that, if $A$, $B$, $C$, $D$ denote the sets of marked points behind
  the 4-valent vertex, the given set $S$ can be written as a union of some of
  the sets $A$, $B$, $C$, $D$. We have to show that up to a global multiple
  there is only one weight distribution on these cones giving rise to a
  tropical cycle. By relabeling the subsets (and possibly replacing $S$ by
  $S^c$) we may assume w.l.o.g.\ that $ S=A $ or $ S=A \cup B $.

  By Lemma \ref {lem-divisor-set} we know already that the weights of these
  cones can only depend on the sets $A$, $B$, $C$, $D$; we denote these weights
  by $ w(A,B,C,D) $. So to prove irreducibility it suffices to show that
  \begin {enumerate}
  \item $ w(A,B,C_1,D \cup C_2) = w(A,B,C,D) $ if $ S=A \cup B $, $ C = C_1 \sqcup C_2 $,
and $ C_1,C_2 \neq \emptyset $,
  \item $ w(A \cup B, C_1,C_2, D) = -w(A,B,C,D) $ if $ S=A \cup B $, $ C= C_1 \sqcup C_2
$, and $ C_1,C_2 \neq \emptyset $,
  \end {enumerate}
  since these two operations can transform any partition of $ \{1,\dots,n\} $
  into four subsets occurring in $ D^S $ into any other.

  To prove this, consider a codimension-1 cone in $ D^S $ having one 5-valent
  vertex, with sets of marked points in $ A,B,C_1,C_2,D $ behind its adjacent
  edges. Since $ S=A \cup B $, only $4$ of the $10$ resolutions of this cone in
  the codimension-1 skeleton of $ \M_{0,n} $ occur in $ D^S $:
  \begin {center} \input {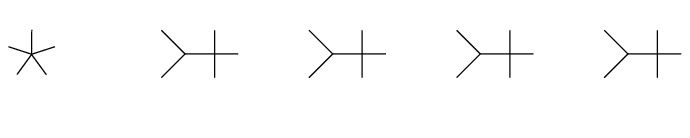} \end {center}
  Now we map these cones to $ \M_{0,5} $ by the forgetful map that sends the
  marked ends in $ A,B,C_1,C_2,D $ to $ x_1,x_2,x_3,x_4,x_5 $, respectively.
  The $4$ cones occurring above are then exactly mapped to the four cones in
  the vital divisor $ D^{\{1,2\}} $ in $ \M_{0,5} $ of Example \ref{ex-vital-5}.
  As this divisor is irreducible (see Example \ref {m05}), we
  conclude that there is up to a global multiple only one choice of weights for
  the four cones above. By the explicit computation of Example \ref{ex-vital-5}
  these weights satisfy $ -w(A \cup B,C_1,C_2,D) = w(A,B,C_1,D
  \cup C_2) = w(A,B,C,D) $, proving both (a) and (b).
\end{proof}

\begin {remark} \label {rem-behavior}
  The proofs of Propositions \ref{pir} and \ref{prop-vital} suggest that, once
  we know that a divisor in $\M_{0,n}$ is locally irreducible, we might be able
  to conclude by some ``connectivity'' argument that the divisor is also
  (globally) irreducible. Unfortunately, the proof of this connectivity is very
  complicated for increasing $n$ since the number of possible cones contained
  in the support of such a divisor increases exponentially. So far, we are just
  able to show the equivalence of local and global irreducibility of a divisor
  in the case $n=6$ (for $n=5$ this is trivial since there is only one
  balancing condition):
\end {remark}

\begin {proposition} \label {prop-m06}
  In $ \M_{0,6} $ every locally irreducible divisor is also globally
  irreducible.
\end {proposition}

\begin{proof}
  Let $D$ be a locally irreducible divisor. We have to show that its weight
  function is unique up to a global multiple in the sense of 
  Lemma \ref{criterion}.

  To prove this, let us introduce the following notation. A codimension-$2$
  cone $\tau$ of $D$ will be called \textit {relevant} if it corresponds to
  curves with one $5$-valent vertex (rather than to curves with two $4$-valent
  vertices). Note that $D$ is determined uniquely by the weights of all cells
  adjacent to relevant cones (hence the name) since every codimension-$1$ cell
  is adjacent to such a relevant codimension-$2$ cone. Moreover, in $ \M_{0,6}
  $ all relevant cones correspond to curves having two markings $i$ and $j$ on
  one side of the bounded edge and four on the other. We will use the notation
  $ij$ for this cone.

  As $D$ is locally irreducible we know that the weight combination of the
  codimension-$1$ cells around a relevant cone $ \tau $ must be a multiple of
  the irreducible types of divisors in $ \M_{0,5} $ listed in Example \ref{m05}.
  If the weights around $ \tau $ are not all $0$ we will call the choice
  of type the \textit {resolution behavior} at $ \tau $, and the corresponding
  multiple the \textit {weight} of $ \tau $. However, in $\M_{0,6}$, where each
  curve has $6$ markings, we have to be careful how we split up the six
  markings over the subtrees behind the five bounded edges adjacent to the
  $5$-valent vertex. This means that all types of Example \ref {m05} correspond
  to several resolution behaviors. The following list shows all possible
  resolution behaviors at the cone $12$ (up to permutation of the labels of the
  markings $ 3,4,5,6 $). In this list, we understand by \textit{new bounded
  edge} the bounded edge appearing when resolving the $5$-valent vertex.
  \begin {itemize}
  \item $\psi_3$, where we require the marking $3$ to be at the $4$-valent
    vertex in the resolving curves;
  \item $\psi_{12}$ (not really a Psi-class!) corresponding to resolving the
    $5$-valent vertex where the subtree consisting of the markings $1$ and $2$
    is required to stay at the $4$-valent vertex;
  \item $V^{\lk 5,6 \rk}$, where either the marking $5$ or $6$ should lie at
    the $4$-valent vertex (but not both);
  \item $V^{\lk 12,3 \rk}$, where the marking $3$ or the subtree consisting of
    the markings $1$ and $2$ should be at the $4$-valent vertex, but not both;
  \item $\psi_3^\natural$, where the marking $3$ is not allowed to lie at the
    $4$-valent vertex;
  \item  $\psi_{12}^\natural$, where the subtree consisting of the markings $1$
    and $2$ is not allowed to lie at the $4$-valent vertex;
  \item  $D^{\lk 5, 6\rk}$, where the markings $5$ and $6$ lie both on the side
    of the new bounded edge without $4$-valent vertex, or both on the side
    of the new bounded edge with the $4$-valent vertex;
  \item $D^{\lk 12, 3\rk}$, where the marking $3$ and the subtree consisting of
    the markings $1$ and $2$ lie both on the side of the new bounded edge
    without $4$-valent vertex, or both on the side of the new bounded edge with
    the $4$-valent vertex.
  \end {itemize}
  For some resolution behaviors, we show the corresponding adjacent
  codimension-$1$ cells occurring in $D$ with non-zero weight in the picture
  below, where unmarked ends are understood to be labeled with the remaining
  markings in any possible way.

  \begin {center} \input {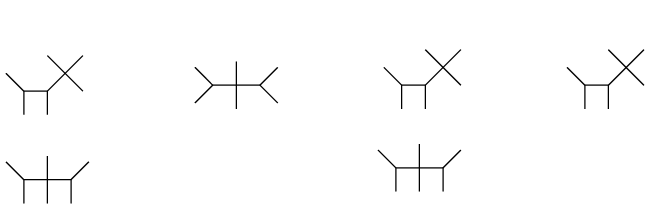} \end {center}

  The main idea of the proof is now the following: pick a relevant
  codimension-$2$ cone, say $12$, such that the weights of $D$ around it are
  not all $0$. If we prescribe the resolution behavior and the weight at $12$
  this will determine the weight (not necessarily the resolution behavior
  however) of all other relevant codimension-$2$ cones adjacent to
  codimension-$1$ cells occurring in the chosen resolution behavior at $12$.
  For example, if we choose the resolution behavior $ \psi_3 $ at $12$ as
  above, the picture below shows the type $12$ (on the left) and an adjacent
  codimension-$1$ cell with non-zero coefficient in $ \psi_3 $ in the middle.
  As this cell is also adjacent to the relevant codimension-$2$ cone $56$ (on
  the right), this means that the weight of $56$ is fixed by that of $12$.

  \begin {center} \input {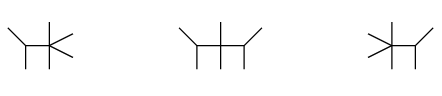} \end {center}

  This argument can be extended slightly by using in addition the result of
  Lemma \ref {lem-divisor-set} that the weights of codimension-$1$ cells in
  divisors depend only on the distribution of the marked points behind the
  $4$-valent vertex. For example, the following picture shows that the weight
  of $14$ is also fixed by the weight (and the resolution behavior) of $12$.

  \begin {center} \input {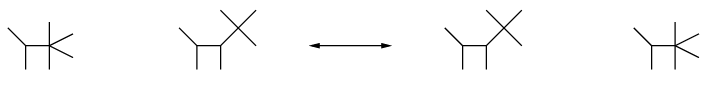} \end {center}

  Using these arguments, the following list shows, when fixing a certain
  resolution behavior at $12$, the weights of which other relevant
  codimension-$2$ cones are fixed by the weight of $12$. The column on the
  right shows the number of relevant codimension-$2$ cones whose weight is
  fixed by that of $12$.
    \[ \begin{array}[c]{|l|l|r|} \hline
         \psi_3 & 12, 14,15,16,24,25,26,45,46,56&10\\ 
         \psi_{12} & 12, 34,35,36,45,46,56&7\\ 
         V^{\lk 5,6 \rk} & 12,15,16,25,26,35,36,45,46&9\\ 
         V^{\lk 12,3 \rk} & 12,14,15,16,24,25,26,34,35,36&10\\ 
         \psi_3^\natural & 12,13,23,34,35,36&6\\ 
         \psi_{12}^\natural & 12,13,14,15,16,23,24,25,26&9\\ 
         D^{\lk 5,6 \rk} & 12, 13, 14,23,24,34,56&7\\ 
         D^{\lk 12,3 \rk} & 12,13,23,45,46,56&6\\ \hline 
       \end{array} \]
  Note that there are $ \binom 62=15 $ relevant codimension-$2$ types in total.
  So if we assume that $D$ is reducible, this means that a resolution behavior of 
  $\psi_3 $ or $ V^{12,3} $ as in the table above is impossible since this fixes
  $10$ weights already, leaving only $5$ of the weights of relevant cones
  undetermined --- which is not enough for another independent weight since
  each weight has to occur in at least $6$ relevant cones by the table above.

  To rule out the other resolution behaviors we cannot use a pure counting argument
  alone, but rather look in detail at the precise cones whose weights are
  determined by a fixed weight in a certain resolution behavior. A straightforward
  but tedious case-by-case analysis using the table above shows that it is not
  possible to fix the resolution behaviors and independent weights at two relevant
  cones without contradictions, i.e.\ such that the sets of relevant cones
  whose weights are determined by the two initial ones are disjoint. This shows
  that $D$ must in fact be irreducible.
\end{proof}

Since the techniques used in Proposition \ref {prop-m06} are limited to $n=6$,
the following question remains open:

\begin{conj}
  Let $n>6$. Is any locally irreducible divisor in $ \M_{0,n} $ also globally
  irreducible?
\end{conj}

  \section {The set of points in special position for rational $n$-marked plane
  tropical curves} \label {l3}

One natural extension of the notion of rational $n$-marked abstract tropical
curve is the rational $n$-marked parametrized tropical curve introduced below,
which we consider to be the tropical analogue of a rational $n$-marked stable
map in classical algebraic geometry. Recall that for a given smooth complex
projective variety $X$, $n \in \N $ and a homology class $\beta \in H_2(X,\Z)$
a \textit{rational $n$-marked stable map (into $X$)} is a collection $(C;
x_1,\ldots,x_n,f)$ of a complete connected nodal curve $C$ of arithmetic genus
$0$, smooth points $x_i \in C$ which are pairwise disjoint, and a morphism $f:C
\to X$ satisfying $f_{∗}([C]) = \beta$ s.t.\ every irreducible component that
is mapped to a point in $X$ contains at least 3 \textit{special points}, i.e.\
markings or intersection points with other components (\textit{stability}). Two
stable $n$-marked maps $(C; x_1,\ldots,x_n, f)$ and $ (\tilde{C};
\tilde{x}_1,\ldots,\tilde{x}_n,\tilde{f}) $ are \textit{isomorphic} if there
exists an isomorphism $\varphi: C \to \tilde{C}$ s.t.\
$\varphi(x_i)=\tilde{x}_i$ and $\tilde{f} \circ \varphi = f$. The space
$\overline{M}_{0,n}(X, \beta)$ is the set of isomorphism classes of rational
$n$-marked stable maps; it carries the structure of a smooth projective variety
\cite{B,FP}.

In both the tropical and the algebraic setup one can then count such curves
through given points, where the number of these points has been chosen so that
one expects a finite answer. We are interested in the description of the set of
points in special position for these counting problems. Remember that in the
classical situation, when the ambient space for the rational stable maps is
$X=\PP^2$, it is known that the points in special position --- in the sense
that there are infinitely many curves passing through them --- form a
subvariety of $(\PP^2)^n$ (see \cite [Corollaire 13.1.5]{EGA} applied to the
(product) evaluation map). In the tropical context a more striking result
holds, namely that these points form a tropical subfan of $(\R^2)^n$ of
codimension one, which we can describe as a push-forward of some divisor in
$\M_{0,n}$. To be more precise, we deal here with two notions of ``points in
special position'' which both as sets arise as a divisor pushed forward by the
evaluation map. We restrict ourselves to the case where the ambient space is
$\R^2$ as our arguments just hold there.

But let us start by recalling the necessary tropical vocabulary. More details
about the following can be found in \cite [section 4]{GKM}.

\begin {definition}[Rational $n$-marked (labeled) parametrized tropical curve
    of degree $\Delta$ in $\R^2$]
   We denote a tuple $(v_1,\ldots,v_m)\in (\Z^2 \setminus \lk 0 \rk)^m$ by
  $\Delta$.

  A \textit {rational $n$-marked (labeled) parametrized tropical curve of
  degree $\Delta$ in $\R^2$} is a pair $(C, h)$ consisting of a rational
  $(n+m)$-marked abstract tropical curve $C=(\Gamma; x_1, \ldots, x_{n+m})$ and
  a continuous map $h:\Gamma \to \R^2$ such that the following holds.
  \begin {itemize}
  \item The map $h$ is integer affine linear on each edge $E$ of $ \Gamma $,
    i.e.\ of the form $ h(t)=a+vt $ for some $ a \in \R^2 $ and $ v \in \Z^2 $.
    If we start parametrizing $E$ at the vertex $ V \in \partial E $ we call
    $v$ the \textit {direction} $ v(E,V) $ of $E$ with respect to $V$.
  \item At each vertex $V$ the \textit {balancing condition}
      \[ \qquad \quad \sum_{E: V \in \partial E} v(E,V) = 0 \]
    holds.
  \item The direction of $ x_i $ is $0$ for all $ i = 1, \ldots,n$ (contracted
    ends).
  \item The direction of $ x_i $ is $ v_{i-n} $ for all $ i=n+1,\ldots,n+m $
    (non-contracted ends).
  \end {itemize}
  Let $V$ be a 3-valent vertex which is not adjacent to a contracted end. Let
  $w_1,w_2,w_3$ be the direction vectors of its adjacent edges. Then we define
  the \textit {multiplicity} $ \mult(V) $ of $V$ to be $ \bet{\det(w_1|w_2)} =
  \bet{\det(w_1|w_3)} = \bet{\det(w_2|w_3)} $.

  In our pictures we will usually only draw the image curve $ h(\Gamma) $
  together with the points $ h(x_1),\dots,h(x_n) $. This image then has $m$
  (labeled) unbounded edges whose directions are contained in $\Delta$. We can
  think of the \textit {degree} $\Delta$ as the tropical equivalent of the
  homology class of an algebraic stable map. Note that the stability condition translates
  to the requirement that $\Gamma$ has at least $3$-valent  vertices.

  In the following we will sometimes abbreviate ``rational $n$-marked (labeled)
  parametrized tropical curve of degree $\Delta$ in $\R^2$'' by ``parametrized
  tropical curve''. 
\end {definition}

\begin{construction} \label{mlab}
  Two parametrized tropical curves $(C,h)$ and
  $(\tilde{C},\tilde{h})$ are
  \textit{isomorphic} if there is an isomorphism $\varphi$ between the underlying abstract
  tropical curves $C$ and $\tilde{C}$ satisfying $\tilde{h} \circ \varphi =h$. The
parameter space of  rational $n$-marked (labeled) parametrized  tropical
  curves of given degree $\Delta$ in $\R^2$ modulo isomorphisms is denoted by
  $\M_{0,n}^{\tn{lab}}(\R^2,\Delta)$. Observe that the directions of the
  bounded edges in a parametrized tropical curve are not fixed by $\Delta$. But
  when fixing a combinatorial type (see Remark \ref{combtypes}) in $\M_{0,n+m}$ for $
m=|\Delta| $ there is
  a unique choice for the directions of the bounded edges for a parametrized
  tropical curve of degree $\Delta$ in $\R^2$ such that the underlying graph
  $\Gamma$ is of this combinatorial type. Hence there is a bijection between
  combinatorial types of the moduli spaces $\M_{0,n}^{\tn{lab}}(\R^2,\Delta)$
  and $\M_{0,n+m}$ \cite[Lemma 4.6]{GKM}. Note that the number of combinatorial
  types in $\M_{0,n}^{\tn{lab}}(\R^2,\Delta)$ is finite too for this reason.

  Coordinates in $\M_{0,n}^{\tn{lab}}(\R^2,\Delta)$ are given by the position
  of a \textit{root vertex} in $\R^2$ of the para\-metrized tropical curve and
  the lengths of the bounded edges (see \cite[Proposition 2.11]{GM2}).
Hence,
  given a parametrized tropical curve in $\M_{0,n}^{\tn{lab}}(\R^2,\Delta)$ we
  can send it to the image under the map $v$ of its underlying graph $\Gamma$
  (see Construction \ref{moduli}). So we get an isomorphism of polyhedral
  complexes \cite[Proposition 4.7]{GKM}
    $$ \M_{0,n}^{\tn{lab}}(\R^2,\Delta) \stackrel{\cong}{\rightarrow}
       v(\M_{0,n+m}) \times \R^2. $$
  In particular, we can consider $\M_{0,n}^{\tn{lab}}(\R^2,\Delta)$ as a
  tropical fan of dimension $(n+m-3)+2$.
\end{construction}

\begin{definition}[More intersection theory and tropical morphisms]
    \label{inter2}
  We extend the definitions of Definition \ref{intersection}.
  \begin {itemize}
  \item[(a)]
    Let $X$ and $Y$ be two $k$-cycles. After possibly adequately refining
    $X$ and $Y$ we can construct a $k$-cycle on $X \cup Y$, called the
    \textit{sum of the cycles $X$ and $Y$}, which is denoted by $X + Y$
    \cite[Construction 2.13]{AR}.
  \item[(b)]
    A \textit {morphism} $f: X \to Y$ \textit{of cycles $X$ and $Y$
    (resp.\ of fans $X$ and $Y$)} is a $\Z$-linear map, i.e.\ a map $\tilde{f}:
    \bet{X} \to \bet{Y}$ induced by a $\Z$-linear map between the underlying
    lattices $\Lambda_X$ and $\Lambda_Y$. 
  \item[(c)]
    Let $f: X \to Y$ be a morphism of a cycle $X$ to an $m$-cycle $Y$, $Z$ an
    $n$-cycle in $X$ where $n \leq m$. Using an appropriate refinement of $X$
    we can assume that the image of each cone in $X$ is a cone of $Y$. We
    define the \textit {push-forward of $Z$ along $f$} by $ f_*(Z) =
    \lk f(\sigma)|~\sigma \in Z\rk$. This polyhedral complex $f_*(Z)$ becomes
    an $n$-cycle in $Y$ by giving the weights
      $$ \qquad \quad
         w_{f_*(Z)} (\sigma') =
         \sum_{\substack{\sigma \in Z \\ f(\sigma)=\sigma'}}
           w_Z(\sigma) \cdot \bet{\Lambda_{\sigma'}/f(\Lambda_\sigma)} $$
    to the $n$-dimensional cones $\sigma' \in Y$ in the image of $Z$ under $f$
    \cite[Proposition 4.6]{AR}, where $\Lambda_\sigma$ is the sublattice of
    $\Lambda_X$ generated by $\sigma$.
  \item[(d)]
    The tropical Psi-classes of Definition \ref {psi} can be written as the
    divisors associated to certain rational functions \cite [Chapter 3]{KM}.
    As such, we can intersect several Psi-classes by consecutively intersecting
    with these rational functions. The product $ \psi_1^{k_1} \cdot \; \cdots
    \; \cdot \psi_n^{k_n} $ is then a cycle whose support consists of all
    curves such that a vertex with the markings $ i_1,\dots,i_m $ has valence
    at least $ k_{i_1}+\cdots+k_{i_m}+3 $ \cite [Chapter 4]{KM}.
  \end {itemize}
\end{definition}

\begin{definition}[Evaluation maps] \label {def-ev}
  For all $i=1,\ldots,n$ define the \textit{$i$-th evaluation map} by
    $$ \ev_i: \M_{0,n}^{\tn{lab}}(\R^2, \Delta) \to \R^2,~(C,h)
       \mapsto h(x_i). $$
  These maps are well-defined and are morphisms of fans
  \cite[Proposition 4.8]{GKM} using Construction \ref {mlab}. The \textit
  {evaluation map}
    $$ \ev:=\prod_{i=1}^n \ev_i: \M_{0,n}^{\tn{lab}}(\R^2,\Delta) \to
       \R^{2n} $$
  is a morphism of fans, too.
\end{definition}

\begin{definition} [Two versions of points in special position]
    \label {def-special}
  Fix $n>0$ and a degree $\Delta$ such that $ n=|\Delta|-1 $ (i.e.\ so that
  the source and target of $\ev$ have the same dimension and we expect a finite
  number of curves of degree $ \Delta $ through $n$ given points). A collection
  $Q=(Q_1,\ldots,Q_n)$ of $n$ points in $\R^2$ is said to be
  \begin {itemize}
  \item \textit {in special position (v1)} if $\ev^{-1}(Q)$ is infinite;
  \item \textit {in special position (v2)} if $\ev^{-1}(Q)$ is infinite or
    intersects the codimension-1 skeleton of $\M_{0,n}^{\tn{lab}}
    (\R^2,\Delta)$.
  \end {itemize}
  Otherwise we say that $Q$ is \textit {in general position} (for (v1) or
  (v2)).
  As $\ev$ is linear on each cone of $\M_{0,n}^{\tn{lab}}(\R^2,\Delta)$, note
  that $\ev^{-1}(Q)$ being infinite is equivalent to saying that the map $\ev$
  is not injective on (at least) one cone of $\M_{0,n}^{\tn{lab}}(\R^2,\Delta)$
  that intersects $ \ev^{-1}(Q) $.
\end{definition}

\begin{remark}
  Version (v2) is typically used for enumerative purposes: fix a point configuration
  $Q$ of points in general position (v2) and consider curves $C$ in $\ev^{-1}(Q)$, which
  contain only $3$-valent vertices. We assign a multiplicity $\mult(C)$ to each
  curve $C$ (see for example \cite{GKM}). It can then be shown that 
  $\sum_{C \in \ev^{-1}(Q)} \mult(C)$ is independent of the choice of $Q$; it
  is an invariant.
\end{remark}
 
\begin{remark} [Welschinger invariants as intersection products]
  Inspired by \cite{MR}, one could hope to write real enumerative invariants,
  for example Welschinger invariants, as intersection products of adequate
  cycles in $\M_{0,n}^{\tn{lab}}(\R^2,\Delta)$. Contrarily to the
  classical situation, the part of $\M_{0,n}^{\tn{lab}}(\R^2,\Delta)$
  corresponding to real $n$-marked rational projective curves has
  the same dimension as $\M_{0,n}^{\tn{lab}}(\R^2,\Delta)$. Hence,
  in order to impose some incidence conditions to these curves
  (like point, line or Psi-conditions) we would have
  to intersect a cycle of codimension $0$ in $\M_{0,n}^{\tn{lab}}(\R^2,\Delta)$
  with pull-backs of rational functions cutting out points, lines,
  or Psi-classes. For real curves, this cycle should be different
  from $\M_{0,n}^{\tn{lab}}(\R^2,\Delta)$, which corresponds to
  the case of complex curves. But Proposition \ref{mir} states that
  such a cycle cannot exist.
\end{remark}

We now want to derive formulas for the locus in $ \R^{2n} $ of points in
special position, for both versions (v1) and (v2).

\begin{remark}[Strings] \label {string}
  A parametrized tropical curve has a \textit{string} if the  underlying graph $\Gamma$
  contains a subgraph homeomorphic to $\R$ which does not  intersect the closures
  $\overline{x_i}$ of the markings $x_1,\ldots,x_n$. By \cite [Remark 3.7]{GM2} curves
  lying in codimension $0$ of $ \M_{0,n}^{\tn{lab}}(\R^2,\Delta) $ passing through points
  in special position have at least one string. Such a curve can have several strings
  which are not necessarily disjoint. 
\end{remark}

\begin {definition}[Free and fixed edges] \label{def-free-fixed}
  Let $(C,h)$ with $ C=(\Gamma; x_1,\ldots,x_n) $ be a parametrized tropical
  curve in $\R^2$. Let $V$ be a vertex of $C$ and $E$ an adjacent edge. Then
  $E$ is called a \textit {free edge} at $V$ if it can be connected in
  $ \Gamma \setminus (V \cup \bigcup_{i=1}^n \overline{x_i})$ to an unmarked
  end. Otherwise we call $E$ a \textit {fixed edge} at $V$.
\end {definition}

\begin{proposition}[Points in special position (v1)] \label{v1}
  Let $n=\bet{\Delta}-1$ and assume $n>1$. Then the set of points in special
  position (v1) for curves in $\M_{0,n}^{\tn{lab}}(\R^2,\Delta)$ equals
  the support of the cycle $\ev_*(\psi_1 +\dots + \psi_n)$.
\end{proposition}

\begin{proof}
  We have to show two inclusions.

  In order to prove that the set of points in special position is contained in
  the support of $ \ev_*(\psi_1+\dots + \psi_n) $ we consider curves lying in
  codimension $0$ of $ \M_{0,n}^{\tn{lab}}(\R^2,\Delta) $ passing through
  points in special position. By Remark \ref{string} these curves have at least
  one string. The idea of the proof is that moving such a string yields a curve
  passing through the same point configuration, but lying in a codimension-one
  cone of $ \M_{0,n}^{\tn{lab}}(\R^2,\Delta) $ and having one $4$-valent vertex
  where one of the adjacent edges is a marking $x_1,\ldots,x_n$. (Remember that
  a marking $x_i$ is an edge of the underlying abstract tropical curve that is
  mapped to a point in the parametrized curve, it is a contracted end.
  So a marking $x_i$ adjacent to a vertex is depicted by a point on the vertex
  below.) This cone is often not unique as we might move the string in
  different directions. Moreover, the resulting curve does not necessarily lie
  in the boundary of the original cone: it might happen that the string first
  runs into a $4$-valent vertex such that no marking $x_i$ is adjacent to it.
  When resolving this vertex a new string appears which can be moved again. The
  following sequence of pictures shows the idea of these movements; there is
  one string which is drawn in dashed lines. In the second picture there is a
  4-valent vertex without adjacent marking. The final curve has a marking at a
  4-valent vertex $ x_i $, so it lies in the support of $ \psi_i $, which means
  that our set of points in special position lies in the support of $ \ev_*
  \psi_i $ and thus also of $ \ev_* (\psi_1+\cdots+\psi_n) $.

  \begin {center} \begin{picture}(0,0)%
\includegraphics{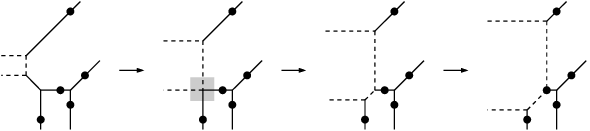}%
\end{picture}%
\setlength{\unitlength}{4144sp}%
\begingroup\makeatletter\ifx\SetFigFont\undefined%
\gdef\SetFigFont#1#2#3#4#5{%
  \reset@font\fontsize{#1}{#2pt}%
  \fontfamily{#3}\fontseries{#4}\fontshape{#5}%
  \selectfont}%
\fi\endgroup%
\begin{picture}(4479,998)(664,-290)
\put(4775, 21){\makebox(0,0)[rb]{\smash{{\SetFigFont{10}{12.0}{\familydefault}{\mddefault}{\updefault}{\color[rgb]{0,0,0}$x_i$}%
}}}}
\end{picture}%
 \end {center}

  To make this argument rigorous we have to give an algorithm how to move a
  string so that it runs into a marking. For this let us first consider a curve
  with a 4-valent vertex on a string without adjacent marking (as in the second
  picture above). The following picture shows the types of 4-valent vertices
  without marking, where the types distinguish which of the adjacent edges are
  parallel.

  \begin {center} \begin{picture}(0,0)%
\includegraphics{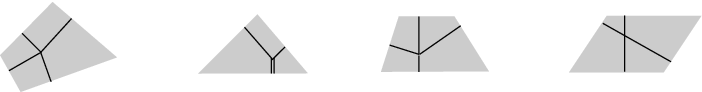}%
\end{picture}%
\setlength{\unitlength}{4144sp}%
\begingroup\makeatletter\ifx\SetFigFont\undefined%
\gdef\SetFigFont#1#2#3#4#5{%
  \reset@font\fontsize{#1}{#2pt}%
  \fontfamily{#3}\fontseries{#4}\fontshape{#5}%
  \selectfont}%
\fi\endgroup%
\begin{picture}(5345,686)(254,-66)
\end{picture}%
 \end {center}

  Note that at most two of the adjacent edges are fixed by the point conditions
  since at least two of them lie on a string. If none of the adjacent edges is fixed,
  the $4$-valent vertex arises from the string movement in curves in codimension $0$
  having at least two strings joining in codimension $1$ at the 4-valent vertex. In
  this case, it is possible to move one of the strings differently in order to obtain a
  4-valent vertex with at least one adjacent fixed edge. Let us assume this in the
  following. Then, considering all possibilities which of the
  edges can be fixed in each of the above types, one can see that in each
  case there is a resolution of the 4-valent vertex such that at least one of
  the fixed adjacent edges becomes shorter. For instance, for the first 4-valent vertex
  from the left in the picture above the table below lists the resolution(s) (A), (B) or
  (C), where at least one of the fixed edges a, b, c or d gets shorter,
  depending on which of the adjacent edges are fixed. Note that, in the case of one fixed
  edge, this edge becomes shorter in each of the resolutions.

  \begin {center} \input {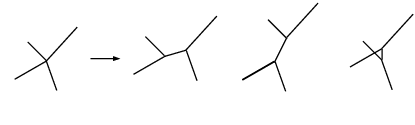} \end {center}
  \begin{center}
    \begin{tabular}{|c||c@{\;\;}|c@{\;\;}|c@{\;\;}|c@{\;\;}|c@{\;\;}|c|} \hline
      edges fixed & a, b & a, c & a, d & b, c & b, d & c, d \\ \hline \hline
      resolutions &(A), (C)&(B)&(A), (B)&(A), (B)&(B), (C)&(A), (C)\\ \hline
    \end{tabular}
  \end{center}

  The algorithm now works as follows. Start with a triple $ (C,S,E) $ where $C$
  is the original curve, $S$ a string on $C$, and $E$ a fixed edge adjacent to
  $S$ (such a choice is possible since there is at least one marking). Denote
  by $ M=M(C,S,E) $ the maximum distance in $C$ from $S$ to a marking behind
  $E$ which can be reached from $S$ without passing other markings. The
  following picture on the left shows such a choice; $S$ is again drawn with
  dashed lines, and $M$ is the length of the two line segments drawn in bold.

  \begin {center} \input {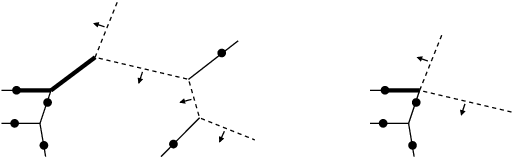} \end {center}

  We now claim that we can always change the curve by moving the string (as
  indicated by the arrows) so that $M$ decreases --- until either $M$ becomes
  zero and thus the string runs into a marking ($x_1$ in the picture above), or
  the string runs into another marking elsewhere earlier (maybe $x_6$ in the
  picture above). The possibility of such a movement is obvious as long as the
  length of $E$ is positive. If $E$ shrinks to a point in the movement (as in
  the picture above on the right) we have a 4-valent vertex at the string with
  at least one fixed adjacent edge, and by our above argument we know that we
  can always continue to move the string so that at least one of the fixed
  adjacent edges becomes shorter. Choosing this edge to be $E$ we can thus
  continue to decrease $M$ (note that by changing $E$ the set of first
  markings behind $E$ is replaced by a smaller one, so this step cannot make
  the maximum $M$ of their distances to $S$ bigger). This completes the
  argument and yields the first inclusion of the proposition.

  For the other direction, we have to show that $\ev_* (\psi_1 + \dots +
  \psi_n)$ contains no points in general position. As $\sum_{i=1}^n \psi_i$ is
  a divisor, we just have to consider curves lying in codimension one of $
  \M_{0,n}^{\tn{lab}}(\R^2,\Delta) $. So consider such a curve $ C =
  (\Gamma;x_1,\dots,x_n) $ in $ \psi_i $, i.e.\ a curve with one $4$-valent
  vertex $V$ with an adjacent marking $x_i$ and only 3-valent vertices
  otherwise. We have to prove that we can deform $C$ to a codimension-0 curve
  that still satisfies the same point conditions.

  This is obvious if $C$ contains a string, so let us assume that this is not
  the case. Note that removing $ \overline {x_i} $ from $ \Gamma $ separates
  $\Gamma$ into $3$ parts, whereas removing each of the other $ n-1 =
  |\Delta|-2 $ causes one more separation. So $ \Gamma \backslash
  \bigcup_{j=1}^n \overline{x_j} $ consists of $ |\Delta|+1 $ connected
  components. As none of these components can have more than one end (otherwise
  we would have a string) we conclude that there is precisely one bounded
  component with no end, whereas all other $|\Delta|$ components contain
  exactly one end. This means that at $V$ (which has $ \psi_i $ and three more
  edges adjacent to it) at least two of the unmarked adjacent edges must be
  connected in $ \Gamma \backslash \bigcup_{j=1}^n \overline{x_j} $ to an
  unbounded edge. We can then resolve $V$ so that these two edges $ E_1 $ and $
  E_2 $ remain together but separate from $ x_i $, forming a string and thus a
  movement of the curve with the positions of the markings fixed.

  \begin {center} \input {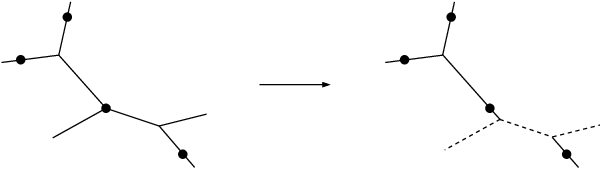} \end {center}

  Hence the points are by definition in special position.
\end{proof}

\begin{corollary} \label {cor-v1}
  In the same situation as above, the set of points in special position
  (v1) for curves in $\M_{0,n}^{\tn{lab}}(\R^2,\Delta)$ can be given the
  structure of a tropical subfan in $\R^{2n}$ of codimension one.
\end{corollary}

\begin{proof}
  Observe that $\ev$ is a morphism of fans of dimension $2n$. The claim then
  follows directly from Proposition \ref{v1} and Definition \ref{inter2} (a)
  and (c).
\end{proof}

\begin{remark}[Comparison to the classical situation]
  As stated in the introduction of this chapter, the set of $n$ points in
  special position (v1) in the corresponding classical situation is a
  subvariety of $(\PP^2)^n$. But in contrast to Corollary \ref {cor-v1} it is
  not necessarily a subvariety of codimension one: consider for instance conics
  in $\PP^2$ through 5 points. There are infinitely many conics through these
  points if and only if two of them coincide or four of them lie on a line ---
  and this forms a subvariety of codimension $2$ in $ (\PP^2)^5 $.

  The reason for the bigger dimension on the tropical side is that there are
  infinitely many liftings in the sense of \cite{JMM} to $\PP^2$ of the points
  in $\R^2$. When tropicalizing, the algebraic curves through each such
  configuration in $(\PP^2)^n$ give rise to tropical curves passing through the
  given configuration of points in $(\R^2)^n$. As a consequence, the number of
  such tropical curves through the given points can be infinite although the
  number of algebraic curves through any lifting of them is not, i.e.\ the
  point configuration in $ (\R^2)^n $ can be in special position (v1) although
  their liftings are not in special position classically.
  
  On the other hand, returning to the case of conics, the locus of points where
  we find reducible curves through them (which roughly corresponds to (v2)) is
  the image of the locus of reducible curves in $ \overline\M_{0,5} (\PP^2,2
  \ell)$ under the evaluation map, where $\ell$ is the homology class of a
  line. Here, both this locus in $ \overline\M_{0,5} (\PP^2,2 \ell)$ and its
  image in $ (\PP^2)^5 $ has codimension 1, the latter being the space of all
  points where three of them lie on a line.
\end{remark}

\begin {remark}[Generalization to curves with Psi-class conditions]
    \label {rem-v1}
  One can generalize the result of Proposition \ref {v1} as follows to the case
  of counting curves satisfying Psi-class conditions (i.e\ to tropical
  descendant Gromov-Witten invariants) as in \cite {MR}: fix $ n>0 $, a degree
  $ \Delta $, and $ k_1,\dots,k_n \ge 0 $ such that $ n =
  |\Delta|-1-k_1-\cdots-k_n $. If we then modify Definitions \ref {def-ev} and
  \ref {def-special} so that they use the moduli space $ \psi_1^{k_1} \cdot \;
  \cdots \; \cdot \psi_n^{k_n} \cdot \M_{0,n}^{\tn{lab}}(\R^2,\Delta) $ instead
  of $ \M_{0,n}^{\tn{lab}}(\R^2,\Delta) $ we count curves through given points
  in $ \R^2 $ with additional Psi-class conditions, i.e.\ such that the valence
  of the vertex with attached marking $ x_i $ is (at least) $ k_i + 3 $ for all
  $i$ \cite [Theorem 4.1]{KM}. The proof of Proposition \ref {v1} can then
  easily be adapted to show that the set of points in special position (v1)
  equals the support of the cycle $ \ev_* (\psi_1^{k_1} \cdot \; \cdots \;
  \cdot \psi_n^{k_n} \cdot (\psi_1 + \cdots + \psi_n)) $. In fact, the first
  direction in the proof of the proposition remains unchanged since it is still
  true that curves passing through points in special position contain a string.
  In the second direction the curves in question will still contain exactly one
  bounded region in $ \Gamma \backslash \bigcup_{j=1}^n \overline{x_j} $; it
  follows that the required movement of the curve is still possible, now
  resolving a $ (k_i+4) $-valent vertex to a $ (k_i+3) $-valent (containing the
  marking) and a 3-valent vertex.
\end {remark}

\begin {proposition}[Points in special position (v2)] \label{v2}
  Fix $\Delta$ and let $n$ equal $\bet{\Delta}-1$. Then the set
  of points in special position (v2) for curves in $ \M_{0,n}^{\tn{lab}}
  (\R^2,\Delta) $ equals the support of the push-forward $ \ev_*
  (\M_{0,n}^{\tn{lab}}(\R^2,\Delta)_{\tn{sk}}^{(2n-1)}) $ of the
  codimension-one skeleton of $ \M_{0,n}^{\tn{lab}} (\R^2,\Delta) $ (note that
  this codimension-one skeleton is a cycle by Proposition \ref {comod} and
  Construction \ref {mlab}).
\end {proposition}

\begin{proof}
  By definition the support of the cycle $ \ev_* (\M_{0,n}^{\tn{lab}}
  (\R^2,\Delta)_{\tn{sk}}^{(2n-1)}) $ contains only points in special position.
  In the other direction, if the points are in special position we can use the
  same argument as in the proof of Proposition \ref {v1} to show that they lie
  in the push-forward of the codimension-one skeleton of $ \M_{0,n}^{\tn{lab}}
  (\R^2,\Delta) $.
\end{proof}

\begin{corollary} \label {cor-v2}
  In the same situation as above, the set of points in special position (v2)
  for curves in $\M_{0,n}^{\tn{lab}}(\R^2,\Delta)$ can be given the structure
  of a tropical subfan in $ \R^{2n} $ of codimension one.
\end{corollary}

\begin{proof}
  Observe that $\ev$ is a morphism of tropical fans of dimension $2n$. The
  claim thus follows directly from Proposition \ref {v2} and Definition \ref
  {inter2} (c).
\end{proof}

\begin{remark} [Reducible cycles]\label{rem-irr}
  Using the characterization of Proposition \ref{v1}, the set of points in special
  position (v1) for curves in $\M_{0,n}^{\tn{lab}}(\R^2,\Delta)$ cannot be an irreducible
  cycle of $\R^{2n}$ as it is the push-forward of the reducible divisor $\psi_1 + \cdots
  + \psi_n$. Likewise, the set of points in special position (v2) for curves in
  $\M_{0,n}^{\tn{lab}}(\R^2,\Delta)$ is not an irreducible cycle as Example
  \ref{petersen} shows.
\end{remark}

\begin{remark} [Psi-condition] \label {rem-v2}
  As in Remark \ref {rem-v1} we can easily generalize Proposition \ref {v2} to
  the case of curves satisfying a Psi-condition $ \psi_i $ in addition to
  incidence conditions with points; the result is then that the set of points
  in special position is the push-forward by $ \ev $ of the codimension-one
  skeleton of $ \psi_i $ (which is a cycle by Proposition \ref {copsi}).
  However, for more than one Psi-condition there is no similar statement as the
  codimension-one skeleton of a product of Psi-classes does not have a
  canonical choice of weights.
\end {remark}

\begin{example} \label{petersen}
  Consider $\M_{0,2}^{\tn{lab}}(\R^2,\Delta)$ with $\Delta=\lk e_1
  +e_2,-e_1,-e_2\rk$ where $e_i$ is the $i$-th standard basis vector of $\R^2$
  (i.e.\ curves of degree one in the notation of \cite[Definition 4.1]{GKM}).
  Then $ \M_{0,2}^{\tn{lab}}(\R^2,\Delta) = \M_{0,5} \times \R^2 $ by
  Construction \ref {mlab}. The space $\M_{0,5}$ can be represented by the
  Petersen graph, depicted below twice. The \textit{Petersen graph} is a
  nonplanar, connected graph with $10$ vertices denoted here $ij$ with $i,j \in
  \{1,\ldots,5\} $ and $ i<j $, which are connected by $15$ edges s.t.\ a
  vertex $ij$ is only linked to the three vertices $kl$ with $k,l \in
  \{1,\ldots,5\}\setminus \{i,j\}$. In this graph, the two-dimensional cones of
  $ \M_{0,5} $ appear as edges, and the one-dimensional cones as vertices. The
  vertex $ ij $ corresponds to the ray of $ \M_{0,5} $ generated by the vector
  $v(\lk i,j \rk)$.

  $\M_{0,2}^{\tn{lab}}(\R^2,\Delta)$ is particularly interesting as in this
  case the sets of points in special position of the two versions coincide,
  more precisely they even coincide as tropical fans with the weights of
  Corollaries \ref {cor-v1} and \ref {cor-v2}. To see this, we observe that the
  codimension-one skeleton of $\M_{0,5}$ contains the cone $12$ which is not
  contained in $\psi_1+\psi_2$. Furthermore, the cones $34$, $35$, and $45$
  have each weight $2$ in the fan $\psi_1+\psi_2$ but just weight $1$ each in
  the codimension-one skeleton of $\M_{0,5}$. But these cones $12$, $34$, $35$,
  $45$ vanish when pushed forward by $ \ev $, since $ \ev $ is not injective on
  them (in fact their images consist of the configurations of two equal
  points in $ \R^2 $ and thus have codimension 2 in $ \R^2 \times \R^2 $).
  
  For more markings, the notions (v1) and (v2) will in general differ.

  \begin {center} \scalebox{0.8}{\input{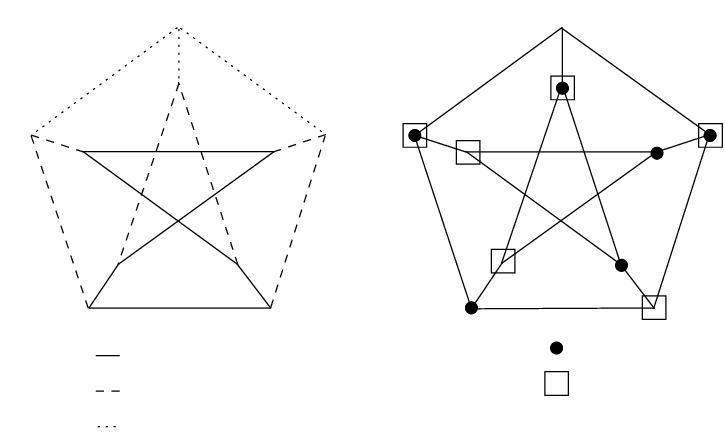}} \end {center} \hfill $\circ$
\end{example}

  To be able to interpret Propositions \ref {v1} and \ref {v2} numerically we now
want to compute the weights of the cells of codimension-1 cycles in $ \R^{2n} $
that are of the form $ \ev_* (Z) $ for a cycle $Z$ in $ \M_{0,n}^{\tn{lab}}
(\R^2,\Delta)$. For this we first need to recall some well-known linear algebra
results on elementary divisors.

\begin{lemdef}[Theorem on elementary divisors] \label{ele}
  Let $M$ be a finitely generated free module over a principal ideal domain
  $R$, and $N \subset M$ a submodule of $M$. Then there exists a basis $ \lk
  u_1,\ldots,u_m \rk$ of $M$, a basis $\lk v_1,\ldots,v_n \rk$ of $N$ and
  $e_1,\ldots,e_n \in R \setminus \lk 0\rk$ such that $v_i=u_i e_i$ for $
  i=1,\dots,n $ and $e_{i+1} \equiv 0 \mod e_i$ for $ i=1,\dots,n-1 $. The
  $e_i$ are called the \textit{elementary divisors of} $N$ and are unique up to
  units in $R$.

  For the quotient module it follows that $M/N \cong R^{m-n} \oplus
  \bigoplus_{i=1}^n R/ R e_i.$ Hence in the case $R=\Z$ the number of elements
  of the torsion part of $ M/N $ is $ \bet{\prod_{i=1}^n e_i} $. In the
  following we will denote this number by $ D(M/N) $.
\end{lemdef}

\begin {proof}
  See for example \cite [Theorem III.7.8]{L}.
\end {proof}

In the rest of the paper we will always use this result for the ring $ R=\Z $.
We then choose the $e_i$ to be positive.

\begin{lemma} \label {lem-gcd}
  In the situation of Lemma \ref{ele} (for $ R=\Z $) the number $ D(M/N) $
  is the greatest common divisor (gcd) of the $ n \times n $
  minors of any matrix $A$ representing the $ \Z $-linear map $ N
  \hookrightarrow M $. 
\end{lemma}
We therefore denote this number by $ D(A) $.
\begin{proof}
  See Remark 3 of Chapter 12.2 on page $6$ of \cite{vdW}. From Lemma \ref
  {ele} it follows that the map $ N \hookrightarrow M $ can be represented by a
  matrix $B$ with the elementary divisors on the diagonal and all other entries
  zero. As $A$ represents the same map there exist matrices $ S \in \tn{GL}
  (m,\Z) $ and $T \in \tn{GL}(n,\Z)$ such that $A=S \, B \, T$. The $n \times
  n$ minors of $A$ are then integer linear functions of the $n \times n$ minors
  of $B$, and vice versa. As the $\gcd$ of the $n \times n$ minors of $B$
  equals the product $ D(M/N) = \prod_{i=1}^n e_i$, this means that, up to
  units in $\Z$, the $\gcd$ of the $n \times n$ minors of $A$ equals $ D(M/N) $
  as well.
\end{proof}

\begin{corollary} \label{3}
  Let $B \in \Z^{(n+1) \times n}$ be a matrix having a $(n-k)\times (n-k)$ part $B_1$ in
  the upper left corner, a block $B_2$ of dimensions $(k+1) \times k$ in the
  lower right corner, an arbitrary (non-quadratic) block in the upper right
  corner, and just zeros in the lower left corner: 
    $$ B= \left( \begin{array}{c|c}
         B_1 & * \\ \hline
         0 & B_2
       \end{array} \right). $$
  Then
    \[ D(B) = |\det (B_1)| \cdot D(B_2). \]
\end{corollary}

\begin{proof}
  To compute the $n \times n$ minors of $B$, we have to erase one row of $B$ and
  look at the determinants of these matrices. If we delete one of the first
  $n-k$ rows, the vectors in the columns of the quadratic part that remain are
  linearly dependent, hence these minors vanish. Deleting the $j$-th row with $
  n-k+1 \le j \le n+1 $, we obtain a matrix that contains the block $B_1$ in
  the upper left corner, a quadratic block $C_j$ of dimensions $ k \times k $
  in the lower right corner and 0 in the lower left part. Hence the determinant
  of such a matrix equals the product of $\det(B_1)$ and $\det(C_j)$. So by
  Lemma \ref {lem-gcd} we get
  \begin{align*}
    D(B)
      &= \gcd \left\{ \det(B_1)\cdot \det(C_j) :\; j = n-k+1,\ldots,n+1
         \right\} \\
      &= |\det(B_1)| \cdot \gcd \left\{ \det(C_j) :\; j=n-k+1,\ldots,n+1
         \right\} \\
      &= |\det(B_1)| \cdot D(B_2),
  \end{align*}
  using the property $ \gcd(ma, mb) = m \, \gcd(a,b) $ for $ m \in \N $ and $
  a,b \in \Z $.
\end{proof}

We will now apply these results to obtain formulas for the weights of
push-forwards of codimension-1 cycles along the evaluation map. For this we
first have to classify those cycles.

\begin {remark}[Codimension-1 types in $ \M_{0,n}^{\tn{lab}}(\R^2,\Delta)
    $] \label {rem-types}
  In the following, a connected component of $ \Gamma \backslash
  \bigcup\limits_{j=1}^n \overline{x_j}$ will be called a \textit {region} of
  the curve. Consider a cell of a codimension-1 cycle in $ \M_{0,n}^{\tn{lab}}
  (\R^2,\Delta) $ on which the evaluation map is injective. It corresponds to a
  combinatorial type of curves having exactly one 4-valent vertex, with all
  other vertices being 3-valent.

  If this 4-valent vertex has an adjacent marking, an argument as in the proof
  of the second part of Proposition \ref {v1} shows that there is exactly one
  region that is bounded (by markings), whereas the others contain
  exactly one end. We will call this type (A); in the picture below the bounded region is
  drawn with dotted lines.

  If the 4-valent vertex has no adjacent marking it lies in a unique region.
  The same argument as above then shows that
  \begin {itemize}
  \item either all regions have exactly one end (type (B) below, with the
    region containing the 4-valent vertex drawn with dotted lines); or
  \item the region with the 4-valent vertex has exactly two ends, there is one
    other bounded region, and all other regions have exactly one end (type (C)
    below, where the bounded region is drawn with dotted lines).
  \end {itemize}
  \begin {center} \input {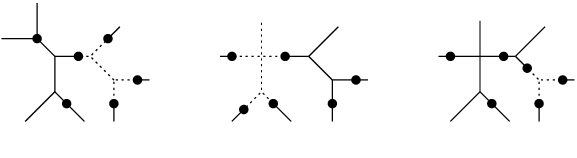} \end {center}
\end {remark}

The weights of the images of those cycles under the evaluation map can be
computed by the following formula.

\begin{lemdef} \label{index}
  Let $Z$ be a cycle of dimension $ 2n-1 $ in $ \M_{0,n}^{\tn{lab}}
  (\R^2,\Delta) $ and $\sigma' \in \ev_*(Z) \subset \R^{2n}$ a cone of the same
  dimension, i.e.\ of codimension 1 in $ \R^{2n} $. Then, with the convention
  and notations from Definition \ref {inter2} (c), the weight of $ \sigma' $ in
  the cycle $ \ev_* (Z) $ is
    $$ w_{\ev_*(Z)}(\sigma')
       = \sum_{\substack{\sigma \in Z \\ \ev(\sigma)=\sigma'}}
         w_Z(\sigma) \cdot \mult(\sigma) $$
  where the multiplicity of a cone $ \sigma $ of dimension $ 2n-1 $ in $
  \M_{0,n}^{\tn{lab}}(\R^2,\Delta) $ is defined as
    \[ \mult (\sigma) := D(\Z^{2n}/\ev(\Lambda_\sigma)). \]
  If $C$ is any curve corresponding to a point in $ \sigma $, we will write the
  multiplicity $ \mult (\sigma) $ also as $ \mult (C) $.
\end{lemdef}

\begin{proof}
  By Definition \ref {inter2} (c) we have to sum over the numbers $ w_Z(\sigma)
  \cdot \bet{\Lambda_{\sigma'}/\ev(\Lambda_\sigma)} $ for all $ \sigma
  \in Z $ with $ \ev(\sigma)=\sigma' $. But as $ \Lambda_{\sigma'} $ is by
  definition a saturated lattice in $ \Z^{2n} \subset \R^{2n} $, we have $
  \Z^{2n} = \Lambda_{\sigma'} \oplus \Z $, and thus the torsion parts of
  $ \Lambda_{\sigma'}/\ev(\Lambda_\sigma) $ and $ \Z^{2n} / \ev(\Lambda_\sigma)
  $ agree.
\end{proof}

\begin {remark} \label {length}
  By Lemma \ref {ele}, the number $ \mult (\sigma) = D(\Z^{2n}/\ev
  (\Lambda_\sigma)) $ in Lemma \ref {index} can be computed as the $ D(A) $ of
  any $ (2n)\times(2n-1) $ matrix $A$ representing the $ \Z $-linear map $ \ev:
  \Lambda_\sigma \to \Z^{2n} $. By \cite[Remark 5.2]{GKM} and
  \cite[Example 3.3]{GM2} one possibility to set up this matrix is to use the
  lengths of all bounded edges and the position in $ \R^2 $ of a root vertex as
  coordinates for $ \Lambda_\sigma $.
\end {remark}

We will now compute explicitly the multiplicities in Lemma \ref {index} for the
cases of Remark \ref {rem-types} and proceed in several steps.

\begin {lemma}[Splitting off vertex multiplicities] \label {lem-reduce}
  Let $C$ be a curve in a codimension-1 cone of $
  \M_{0,n}^{\tn{lab}}(\R^2,\Delta) $ as in Remark \ref {rem-types}. Assume that
  there is an edge of $C$ such that splitting this edge yields two parts $ C_1
  $ and $ C_2 $ of $C$, where $C_2 $ contains only regions with exactly one end and having
  only 3-valent vertices.
  Then
    \[ \mult (C) = \mult (C_1) \cdot \prod_{V \in C_2} \mult (V), \]
  where the product runs over all vertices in $ C_2 $ with no adjacent marking,
  and the multiplicity $ \mult (V) $ of such a vertex $V$ is defined as usual
  as the absolute value of the determinant of two of the adjacent direction
  vectors \cite [Definition 3.5]{GM2}. The following picture shows an example.

  \begin {center} \input {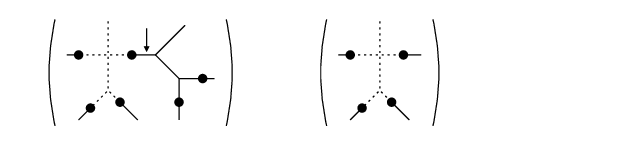} \end {center}
\end {lemma}

\begin {proof}
  Let $k$ be the number of unbounded ends of $ C_2 $. Then $ C_2 $ has $ k-1 $
  markings and $ 2k-2 $ bounded edges (including the split edge). Choosing the
  root vertex to be in $ C_1 $ (see Remark \ref {length}), only the $ 2k-2 $
  coordinates of the $ k-1 $ markings in $ C_2 $ depend on the $ 2k-2 $ lengths
  of the bounded edges in $ C_2 $. Hence the matrix for $ \ev $ as in Remark
  \ref {length} has the form as in Corollary \ref {3}, with $ B_1 $ the $
  (2k-2)\times(2k-2) $ block consisting of these coordinates and lengths. As
  the absolute value of the determinant of $ B_1 $ equals the product of the
  vertex multiplicities in $ C_2 $ by \cite [Proposition 3.8]{GM2} and $ B_2 $
  is precisely the matrix for the evaluation map on $ C_1 $, the claim follows
  from Corollary \ref {3}.
\end {proof}

\begin {lemma}[Multiplicity of a bounded region] \label {lem-bounded}
  Let $C$ be a $3$-valent curve that has exactly one bounded region, with all
  other regions being single ends (such a curve occurs for instance as a part of the
  types (A) and (C) in Remark \ref {rem-types}). Then the multiplicity of $C$
  is
    \[ \mult (C) = \gcd \{ w(E) :\; \text {$E$ end in $C$} \}
         \cdot \prod_{V \in C} \mult (V), \]
  where the weight $ w(E) $ of the edge $E$ is the gcd of the two coordinates
  of the direction vector of $E$, and the product is taken over all vertices of
  $C$ with no adjacent marking. The following picture shows an example.

  \begin {center} \input {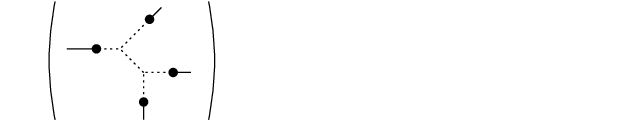} \end {center}
\end {lemma}

\begin {proof}
  We set up the matrix $A$ for the evaluation map as in Remark \ref {length},
  with the root vertex within the bounded region. To compute the minors of $A$
  as required by Lemma \ref {lem-gcd} we have to erase one of its rows. Note
  that the rows of $A$ correspond to the coordinates of the markings in $ \R^2
  $. So let us assume that we erase the row for the $i$-th coordinate of the
  marking $ x_j $ for some $ i=1,2 $ and $ j=1,\dots,n $. Note that the length
  of the bounded edge $ E_j $ adjacent to $ x_j $ is needed only for the
  coordinates of $ x_j $ in $ \R^2 $, and so in the remaining matrix the column
  corresponding to $ E_j $ has at most one non-zero entry, namely for the other
  coordinate of $ x_j $ in $ \R^2 $. Laplace expansion of the determinant
  w.r.t.\ the $ E_j $ column thus simply gives the product of this coordinate
  of $ x_j $ in $ \R^2 $ and the determinant of the evaluation matrix for the
  curve where the marking $ x_j $ is deleted (and thus $ E_j $ becomes an
  unbounded end). But this determinant just equals the product of all vertex
  multiplicities by \cite [Proposition 3.8]{GM2}.

  Altogether we see that $ D(A) $ is the product of all vertex multiplicities
  times the gcd of both coordinates of all markings in $ \R^2 $, as we have
  claimed.
\end {proof}

\begin {corollary}[Multiplicity of the types (A) and (C)] \label {cor-ac}
  Let $C$ be a curve in codimension 1 as in Remark \ref {rem-types}. If $C$ is
  of type (A) with bounded region $ C_b $ then its multiplicity is
    \[ \mult (C) = \gcd \{ w(E):\; \text {$E$ edge in $ C_b $ with adjacent
       marking} \} \cdot \prod_{V \in C} \mult(V). \]
  If it is of type (C) with bounded region $ C_b $ then its multiplicity is
    \[ \mult (C) = \gcd \{ w(E):\; \text {$E$ edge in $ C_b $ with adjacent
       marking} \} \cdot |\det (v,v')| \cdot \prod_{V \in C} \mult(V) \]
  where $v$ and $v'$ are the directions of the two fixed adjacent edges at the
  4-valent vertex (i.e.\ the ones that do not connect to an end within their
  region when coming from the 4-valent vertex). In both formulas, the product
  is taken over all 3-valent vertices without adjacent marking in $C$.
\end {corollary}

\begin {proof}
  Let $V$ be the 4-valent vertex of $C$. If $C$ is of type (A) we can first use
  Lemma \ref {lem-reduce} to split off all vertices behind the two unmarked
  edges adjacent to $V$ that do not lead to the bounded region. This way we get
  the multiplicities of all split-off vertices as a factor, and are left with
  a curve where two of the unmarked edges adjacent to $V$ are solitary ends (as
  it is already the case in the example picture in Remark \ref {rem-types}).
  Now the evaluation matrix of this curve is precisely the same as for the
  curve where these two ends with direction vectors $ v_1 $ and $ v_2 $ are
  replaced by one end with direction $ v_1+v_2 $. Now in the remaining curve we
  can continue to split off all vertices that lie outside of the (closure of
  the) bounded region. This way we are left with a curve whose multiplicity has
  been computed in Lemma \ref {lem-bounded}. Altogether, we get the result
  stated in the corollary.

  If $C$ is of type (C) the procedure is very similar. We first split off all
  vertices behind the two free edges adjacent to $V$ and replace the resulting
  two solitary ends at $V$ by one. This makes the 4-valent vertex $V$ into a
  new 3-valent one for which two adjacent direction vectors are $v$ and $v'$.
  As above, we continue to split off all vertices that are outside of the
  bounded region (one of which will be the new one with multiplicity $ |\det
  (v,v')| $), and use Lemma \ref {lem-bounded} to obtain the result.
\end {proof}

\begin {lemma}[Multiplicity of a region with 4-valent vertex and one end]
    \label {lem-4valent}
  Let $C$ be a curve in codimension 1 as in Remark \ref {rem-types} that has a
  region with a 4-valent vertex and one end directly adjacent to it, all other
  vertices being 3-valent and all other regions being single ends. Then the
  multiplicity of $C$ is
    \[ \mult (C) = \gcd \{ w(E) \cdot |\det (v^1_E,v^2_E)| :\;
         \text {$E$ end with marking in $C$} \}
         \cdot \prod_{V \in C} \mult (V), \]
  where the product is taken over all 3-valent vertices of $C$ with no adjacent
  marking, and $ v^1_E $ and $ v^2_E $ denote the direction vectors of the two
  fixed edges adjacent to the 4-valent vertex that do not connect to $E$. The
  following picture shows an example.

  \begin {center} \input {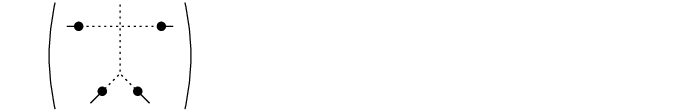} \end {center}
\end {lemma}

\begin {proof}
  The proof is very similar to that of Lemma \ref {lem-bounded}. Let $V'$ be
  the $4$-valent vertex of $C$ and $E$ its unique free end. We set up the
  matrix $A$ of the evaluation map using $V'$ as the root vertex. To compute a
  maximal minor of $A$ we delete the row corresponding to the $i$-th coordinate
  of the marking $ x_j $. Performing a Laplace expansion of the minor w.r.t.\
  the column corresponding to the length of the bounded edge adjacent to $ x_j
  $ we obtain the other coordinate of $ x_j $ in $ \R^2 $ times the determinant
  of the evaluation matrix corresponding to the curve where the marking $ x_j $
  has been deleted. In this new curve we can use the technique of Lemma \ref
  {lem-reduce} to split off all vertices behind the one that lead to $ x_j $.
  We can then replace the two resulting solitary ends at $V'$ ($E$ and the one
  just created by splitting off vertices) at $V'$ by one, leading to a new
  3-valent vertex with multiplicity $ |\det (v^1_E,v^2_E)| $. The resulting
  determinant gives the product of all vertex multiplicities by 
  \cite[Proposition 3.8]{GM2}. Taking the gcd of these expressions for all rows of
  $A$ yields the desired result.
\end {proof}

\begin {corollary}[Multiplicity of the type (B)] \label {cor-b}
  Let $C$ be a curve in codimension 1 as in Remark \ref {rem-types}. If $C$ is
  of type (B), and $C'$ denotes the region with the 4-valent vertex, then its
  multiplicity is
    \[ \mult (C) = \gcd \{ w(E) \cdot |\det (v^1_E,v^2_E)| \} \cdot
       \prod_{V \in C} \mult(V) \]
  where
  \begin {itemize}
  \item the product is taken over all 3-valent vertices without adjacent
    marking in $C$;
  \item the $ \gcd $ is taken over all edges in $C'$ that are adjacent to a
    marking and lie behind one of the three fixed edges adjacent to the
    4-valent vertex; and
  \item $ v^1_E $ and $ v^2_E $ denote as in Lemma \ref {lem-4valent} the
    directions of the two fixed edges adjacent to the 4-valent vertex that do
    not connect to $E$.
  \end {itemize}
\end {corollary}

\begin {proof}
  The proof is similar to that of Corollary \ref {cor-ac}: first we can split
  off all vertices outside of $C'$ and behind the one free end of the 4-valent
  vertex, and then we use Lemma \ref {lem-4valent} for the resulting curve.
\end {proof}

Summarizing, we can now rephrase Lemma \ref{index} as follows.

\begin{corollary} \label{resum}
  Let $Z$ be a cycle of dimension $ 2n-1 $ in $ \M_{0,n}^{\tn{lab}}
  (\R^2,\Delta) $ and $\sigma' \in \ev_*(Z) \subset \R^{2n}$ a cone of the same
  dimension, i.e.\ of codimension 1 in $ \R^{2n} $. Then, with the convention
  and notations from Definition \ref {inter2} (c), the weight of $ \sigma' $ in
  the cycle $ \ev_* (Z) $ is
    $$ w_{\ev_*(Z)}(\sigma')
       = \sum_{\substack{\sigma \in Z \\ \ev(\sigma)=\sigma'}}
         w_Z(\sigma) \cdot \mult(\sigma) $$
  where the multiplicity $ \mult(\sigma) $ of a cone $ \sigma $ of dimension $
  2n-1 $ in $ \M_{0,n}^{\tn{lab}}(\R^2,\Delta) $ is given by Corollaries
  \ref{cor-ac} and \ref{cor-b} depending on their type as in Remark
  \ref{rem-types}.

  In particular, this gives an explicit formula for the cycles in $ \R^{2n} $
  of points in special position (v1) and (v2) as in Propositions \ref {v1} and
  \ref {v2}.
\end{corollary}

\end {document}